\input amstex
\documentstyle{amsppt}
\magnification=\magstep1
\hsize=5.2in
\vsize=6.8in
\input amstex
\documentstyle{amsppt}
\magnification=\magstep1
\hsize=5.2in
\vsize=6.8in

\topmatter
 
\centerline  {\bf RIGIDITY RESULTS FOR WREATH}
\centerline {\bf  PRODUCT II$_1$ FACTORS}

\vskip .1in
\centerline {\rm ADRIAN IOANA}
\address Math Dept  
UCLA, Los Angeles, CA 90095-155505\endaddress
\email adiioana\@math.ucla.edu \endemail
\topmatter
\abstract
We consider II$_1$ factors of the form $M=\overline{\bigotimes}_{G}N\rtimes G$, where either i) $N$ is a non-hyperfinite II$_1$ factor and $G$ is an ICC amenable group or ii) $N$ is a weakly rigid II$_1$ factor and $G$ is ICC group and where $G$ acts on $\overline{\bigotimes}_{G}N$ by Bernoulli shifts.  We prove that isomorphism  of two such factors implies cocycle conjugacy of the corresponding Bernoulli shift actions. In particular, the groups acting are isomorphic. As a consequence, we can distinguish between certain classes of group von Neumann algebras associated to wreath product groups.
\endabstract
\endtopmatter
\document

\head 0. Introduction. \endhead
In recent years, Popa's deformation vs rigidity approach to the classification  of von Neumann algebras has led to  a series of striking results ([Po1,2,3,5,6]). 
Following his terminology, a von Neumann subalgebra $Q$ of a finite von Neumann algebra $N$ is called {\it relatively rigid}
  if any deformation of $id_N$ by subunital, subtracial, completely positive maps $\{\phi_n\}_n$ converges uniformly to $id_N$ on the unit ball of $Q$ ([Po1]).
In the case of algebras $N$ which admit many deformations, Popa's strategy is to use these deformations in combination with the relative rigidity of $Q$.
This in turn gives information about the subalgebra $Q$, which, ideally, reduces the study of the single algebra $N$ to that of the inclusion $Q\subset N$.

A class of von Neumann algebras to which this philosophy has been successfully applied is that of crossed products associated to classical or non-commutative Bernoulli shift actions ([Po2],[Po3],[Po5]).
In general, if $(B,\tau)$ is a finite von Neumann algebra (referred to as the {\it base}) and $G$ is a countable discrete group then the action $\sigma:G\rightarrow\text{Aut}(\overline{\bigotimes}_{G}B)$ given by $\sigma_g(\otimes_{h}x_h)=\otimes_h x_{g^{-1}h},$ for all $g\in G$  is called the $(B,\tau)$-{\it Bernoulli} $G$-{\it action}. 
By analogy with the group case, we  denote by $B\wr G$ the crossed product $(\overline{\bigotimes}_{G}B)\rtimes_{\sigma}G$ and we call it the {\it wreath product of $N$ by $G$}.

Assume now that the {\it base} $(B,\tau)$ is an abelian von Neumann algebra $L^{\infty}(X,\mu)$, or the hyperfinite II$_1$ factor $R$, in which case we call the $(B,\tau)$-Bernoulli actions  {\it classical}, respectively {\it non-commutative}. Popa showed that such actions feature a powerful deformation property called {\it malleability.} 
In the case of a Bernoulli action of $G$ on $M=\overline{\bigotimes}_{G}B$, this amounts to the existence of a continuous action $\theta$ of $\Bbb R$ on $\tilde M=M\overline{\otimes}M$ such that $\theta_1(M\otimes 1)=1\otimes M$ and $\theta_t$ commutes with the double action $\tilde\sigma=\sigma\otimes\sigma$, for any $t\in\Bbb R$ ([Po2],[Po5]).
Using the malleability property in connection with deformation/rigidity arguments, Popa proved  
that if $Q\subset M\rtimes_{\sigma}G$ is a relatively rigid von Neumann subalgebra whose normalizer generates a factor, then $Q$ can be unitarily conjugated into $L(G)$ ([Po2]).

Turning to Bernoulli actions with the base $(B,\tau)$ an arbitrary finite von Neumann algebra, one can ask whether Popa's result on detecting relatively rigid subalgebras of the associated wreath products still holds true.
To answer this question we observe that although, in general, such actions are not malleable, they do satisfy a weaker deformation property, which we call {\it weak malleability}. 
More precisely, we let $\tilde M=\overline{\bigotimes}_{G}[B*L(\Bbb Z)]$, $\tilde\sigma$ be the Bernoulli $G$-action on $\tilde M$ and $\theta_t\in\text{Aut}(\tilde M)$ be defined by $\theta_t=\otimes_{g\in G}(\text{Ad}(u_t))$, where $u_t=\exp(t\log(u))$  for some Haar unitary $u\in L(\Bbb Z)$ and for any $t\in \Bbb R$. 
Using these deformations and relying heavily on Popa's intertwining bimodule techniques ([Po2]) we prove the following dichotomy (Theorem 0.1.) for relatively rigid subalgebras of wreath product von Neumann algebras. Note that $\theta_t$ commutes with $\tilde\sigma$ for any $t\in\Bbb R$
 and that $\theta_1(M)\perp M$. However, since $\theta_1(M)$ and $M$ do not "boundedly" generate $\tilde M$, the proof is more involved than in the malleable case, treated in [Po2] and [Po5],  being closer in spirit to the proof of 6.1. in [Po5] and of  the main result in [IPP].
 
\proclaim {Theorem 0.1} Let $(B,\tau)$ be a finite von Neumann algebra and $G$ be a countable discrete group.
Denote $M=B\wr G$ and suppose that $Q\subset M$ is a relatively rigid diffuse von Neumann subalgebra. If the normalizer of $Q$ generates a factor $P$, then either we can find a unitary $u\in M$ such that  $uQu^*\subset uPu^*\subset L(G)$ or a corner of $Q$ embedds into $\overline{\bigotimes}_{F}B$, for a finite set $F\subset G$. Moreover, in the second case, if we assume that $B$ is a factor, then we can find $u\in\Cal U(M)$ such that  $uQu^*\subset\overline{\bigotimes}_{F}B$ and $uPu^*\subset\sum_{g\in K}(\overline{\bigotimes}_{G}B)u_g$ for some finite set $K\subset G$.
\endproclaim

Rigidity properties for wreath product groups and von Neumann algebras have also been studied by Martin and Vallette ([MV]) and Hayashi.

Recall that a II$_1$ factor $N$ is called {\it weakly rigid} if it has a diffuse, relatively rigid, regular von Neumann subalgebra. Natural examples of weakly rigid factors are provided by group von Neumann algebras, $L(H)$, associated to ICC  {\it weakly rigid} groups $H$, i.e. groups which have an infinite, normal subgroups $H_0$ such that the pair $(H,H_0)$ has relative property (T) of Kazhdan-Margulis. 

\proclaim {Corollary 0.2} Let $N,P$ be two weakly rigid II$_1$ factors and $G,H$ be two countable discrete groups.
If $\theta:N\wr G\rightarrow (P\wr H)^t$ is an isomorphism for some $t>0$, then we can find $u\in\Cal U((P\wr H)^t)$ such that $\theta(\overline{\bigotimes}_{G}N)=u(\overline{\bigotimes}_{H}P)^tu^*$. In particular, $G\simeq H$.
\endproclaim
On the other hand, if we assume the group $G$ to be weakly rigid, then using methods from [Po5] we can  show that any generalized 1-cocycle  associated to any Bernoulli $G$-action is cohomologous to a character of $G$. In turn, combining this result with Corollary 0.2. we can provide a  new class of II$_1$ factors with trivial fundamental group. The existence of such factors has been an open problem for a long time, until settled in the breakthrough work of Popa ([Po1]).

\proclaim{Corollary 0.3} Let $N$ be a weakly rigid II$_1$ factor  and $G$ be  a discrete group having an infinite relatively rigid subgroup. Then $\Cal F(N\wr G)=\{1\}$.
\endproclaim

In the second part of this paper we exploit a rigidity phenomenon coming from central sequences to prove the statement of Corollary 0.2.  for wreath products $N\wr G$ with $G$ amenable and $N$ a non-hyperfinite II$_1$ factor.  A first observation is  that by Connes's characterization of injectivity ([Co1]), $N$ is a non-hyperfinite II$_1$ factor iff all central sequences in $N\wr G$ assimptotically lie in  $\overline{\bigotimes}_{G}N$. On the other hand, since $G$ is amenable, we can use Ornstein-Weiss' Rokhlin lemma ([OW]) to construct many central sequences in $N\wr G$.
\proclaim{Theorem 0.4} Let $N$ be a II$_1$ factor, let $G$ be an infinite discrete group and denote $M=N\wr G$. Let $\omega$ be a free ultrafilter on $\Bbb N$.

$(i)$ If $N$ is non-hyperfinite, then $N'\cap M^{\omega}\subset (\overline{\bigotimes}_{G}N)^{\omega}$.

$(ii)$ Assume that $G$ is amenable and that $Q\subset M$ is a regular von Neumann sublalgebra such that $M'\cap M^{\omega}\subset Q^{\omega}$. Then a corner of $\overline{\bigotimes}_{G}N$ embedds into $Q$.
\endproclaim 

\proclaim{Corollary 0.5} Let $N,P$ be two non-hyperfinite II$_1$ factors and $G,H$ be two countable discrete groups, one of them  amenable.
 If $\theta:N\wr G\rightarrow P\wr H$ is an isomorphism, then we can find $u\in\Cal U(P\wr H)$ such that $\theta(\overline{\bigotimes}_{G}N)=\overline{\bigotimes}_{H}P$. In particular, $G\simeq H$. 
\endproclaim

This paper is organized as follows: in section 1 we review Popa's techniques for conjugating subalgebras of a finite von Neumann algebra and we derive several consequences. 
Section 2 is devoted to  defining the notion of {\it weak malleability} and to discussing consistency and examples: we prove that it indeed generalizes Popa's {\it malleability} (Proposition 2.2.) and that arbitrary Bernoulli shift actions are {\it weakly malleable} (Proposition 2.3.).

In section 3, we first prove a general statement (Theorem 3.3.) asserting that rigid subalgebras of a crossed product $N\rtimes_{\sigma}G$ coming from a weakly malleable action $\sigma$ must have a corner which embedds either into the core $N$ or into the group algebra $L(G)$. Furthermore, if $\sigma$ is a Bernoulli shift action, then we can make this statement more precise, thereby proving Theorem 0.4.   The proofs	 of the above mentioned vanishing cohomology result (Corollary 4.3.) and  of Corollaries 0.2. and 0.3. are the subject of section 4.

Section 5 deals with central sequences in wreath product factors and uses a result of Jones ([Jo1]) to conclude part $(i)$ of Theorem 0.4 (Corollary 5.4.).
In the last section we complete the proof of Theorem 0.4. and then we derive Corollary 0.5.
\vskip 0.05in

\head 1. Technical results. \endhead

We begin this section by reviewing Popa's intertwining bimodule techniques.  

\proclaim{ Theorem 1.1 ([Po2])} Let $(M,\tau)$ be a finite von Neumann algebra and let $Q,B\subset M$  two von Neumann subalgebras. Consider the following conditions: 

(i) There exists $a\in Q'\cap <M,e_B>,a\geq 0,a\not=0,Tr(a)\leq\infty.$ 

(ii) There exist non-zero projections $p\in\Cal P(B),q\in\Cal P(Q)$, an unital homomorphism $\psi$ of $qQq$  into $pBp$ and non-zero partial isometry $v\in M$ such that $vv^*\in (qQq)'\cap qMq, v^*v\in\psi(qQq)'\cap pMp$ and $xv=v\psi(x),\forall x\in qQq.$ 

(iii) For all $a_1,a_2,..,a_n\in M$, $\varepsilon>0,$ there exists $u\in\Cal U(Q)$ such that $||E_B(a_iua_j^*)||_2\leq\varepsilon,\forall i,j=1,.,n.$

(iii)' For any group $\Cal U \subset \Cal U(Q)$ such that $\Cal U''=Q$ and  for all $a_1,a_2,..,a_n\in M$, $\varepsilon>0$, there exists $u\in\Cal U$ such that $||E_B(a_iua_j^*)||_2\leq\varepsilon,\forall i,j=1,.,n.$

Then (i)$\Longleftrightarrow (ii)\Longleftrightarrow$ non(iii)$\Longleftrightarrow$ non(iii)'. 
\endproclaim
\vskip 0.05in
If condition $(ii)$ in the above theorem holds true we  say that "a corner of $Q$ embedds into $B$ inside $M$".
For most applications we will use Popa's theorem to prove $(ii)$ and to do this we will reason by contradiction, assuming that $(iii)$(or $(iii)'$) holds true.
Next, we note a few corollaries.

\proclaim{ Corollary 1.2} Let $(M,\tau)$ be a finite von Neumann algebra and let $Q_1,
Q_2,B\subset M$ be von Neumann subalgebras such that $Q_1$ is generated by unitaries which normalize $Q_2$ (e.g. if $[Q_1,Q_2]=0$). Let $\varepsilon_1,\varepsilon_2>0$, with $\varepsilon_1+\varepsilon_2<1$ and assume that either

(i) $Q_i\subset_{\varepsilon_i}B,\forall i\in\{1,2\}$ or

(ii) $B\subset M$ is regular (i.e. $\Cal N_{M}(B)''=M$)  and there exist finite dimensional left  Hilbert $B$-modules $\Cal H_i\subset L^2(M)$   such that $Q_i\subset_{\varepsilon_i}\Cal H_i$.

 Then  a corner of $Q=(Q_1\cup Q_2)''$ embedds into $B$ inside $M$.

\endproclaim
{\it Proof.} Let $\Cal U=\{u_1u_2|u_i\in\Cal U(Q_i),i=1,2, u_1Q_2u_1^*=Q_2\}\subset \Cal U(Q)$. Then by the assumptions that we made $\Cal U$ is a group and $\Cal U''=Q$.

$(i)$ Since  $Q_i\subset_{\varepsilon_i}B,i=1,2$, we deduce that $$||u_1u_2-E_B(u_1)E_B(u_2)||_2\leq ||u_1(u_2-E_B(u_2))||_2+||(u_1-E_B(u_1))E_B(u_2)||_2$$ $$\leq \varepsilon_1+\varepsilon_2,\forall u=u_1u_2\in\Cal U,$$
hence $$||E_B(u_1u_2)||_2\geq 1-\varepsilon_1-\varepsilon_2=\delta>0,\forall u=u_1u_2\in\Cal U.$$
Thus, condition $(iii)'$ in the previous theorem fails, hence $(ii)$ holds true.

$(ii)$ Let $\omega$ be a free ultrafilter on $\Bbb N$ and let $\Cal N\subset M^{\omega}$ be the von Neumann algebra generated by $B^{\omega}$ and $\Cal N_{M}(B)$. The hypothesis then implies that $Q_i^{\omega}\subset_{\varepsilon_i}\Cal N,i=1,2$
 and the same estimate as above gives that $$||E_{\Cal N}(U)||_2\geq \delta>0,\forall  U=(u_n)_n,u_n\in\Cal U.$$

Assuming that no corner of $Q$ can be embedded into $B$ inside $M$, we can find $u_n\in \Cal U,n\geq 1$ such that $$
\lim_{n\rightarrow \infty}||E_B(xu_n)||_2=0,\forall x\in M.$$ If we denote $U=(u_n)_n$, then $U\perp xB^{\omega},\forall x\in M,$
thus $E_{\Cal N}(U)=0$, a contradiction.
\hfill $\blacksquare$
\vskip 0.05in
 
The first part of the following result is a consequence of  the proof of Proposition 12 in [OP]. We refer the reader to this paper for a proof. 
\proclaim{ Corollary 1.3 ([OP])} Let $M$ be a II$_1$ factor and  $N\subset M$  a subfactor. Let $B=N'\cap M$ and assume that $Q\subset M$ is a von Neumann subalgebra such that a corner of it embedds into $B$.

(i) If $\Cal N_{M}(Q)''$ and $B$ are factors, then there exists $t>0$, a decomposition $B\overline{\otimes}N=B^t\overline{\otimes}N^{1/t}$ and  $u\in\Cal U(M)$ such that $uQu^*\subset B^t$.

(ii) If $Q$ and $B$ are abelian and $\Cal N_{M}(Q)''$ is a factor, then there exists $u\in\Cal U(M)$ and a partition of unity  $1=\sum_{i=1}^{n} q_i$ in $N$ such that $uQu^*\subset \bigoplus_{i=1}^{n}Bq_i$.

\endproclaim 
\vskip 0.1in
If $B$ is a finite von Neumann algebra with a faithful normal trace $\tau$ and $X$ is a countable set, then we denote by $\overline{\bigotimes}_{X}B$, the infinite tensor product $\overline{\otimes}_{x\in X}(B)_x$ with its natural trace, $\otimes_{x\in X}(\tau)_x$.
\vskip 0.05in
{\it Definitions 1.4.} Let $(B,\tau)$ be a finite von Neumann algebra and $G$ be an infinite discrete group acting on a countable set $X$. The action $\sigma:G\rightarrow \text{Aut}(\overline{\bigotimes}_{X}B)$ given by $\sigma(g)(\otimes(a_x)_x)=\otimes(a_{g^{-1}x})_x,$ for all $g\in G$ and $ a=(a_x)_x\in\overline{\bigotimes}_{X}B$ , is called the $(B,\tau)$-{\it Bernoulli} $(G\curvearrowright K)$-{\it action} ([Po5]).  
The {\it wreath product } $B\wr_{X}G$ is then defined to be the crossed product von Neumann algebra  $(\overline{\bigotimes}_{X}B)\rtimes_{\sigma}G.$  Note that if $H$ is a countable discrete group, then  $L(H)\wr_{X}G=L(H\wr_{X}G)$, where $H\wr_{X}G$ is the {\it wreath product} group $H^X\rtimes G$.  

In the case $X=G$ (endowed with the left action of $G$ given by mutiplication), we define the {\it wreath product of $B$ by $G$}
to be $B\wr G=(\overline{\bigotimes}_{G}B)\rtimes_{\sigma}G$ and we call $\sigma$  the $(B,\tau)$-{\it Bernoulli} $G$-{\it action}. 

\vskip 0.05in
The next lemma gives control over the normalizers of certain subalgebras of $\overline{\bigotimes}_{X}B$ inside $B\wr_{X}G$. More generally, we consider product actions $\pi=\sigma\times\beta:G\rightarrow \text{Aut}((\overline{\bigotimes}_{X}B)\overline{\otimes}P)$ given by $\pi(g)=\sigma(g)\otimes\beta(g)$, for all $g\in G$, where $\sigma$ is defined as above and $\beta:G\rightarrow\text{Aut}(P,\tau)$ is an action on a finite von Neumann algebra $(P,\tau)$.
\proclaim{Lemma 1.5} Let   $M=[(\overline{\bigotimes}_{X}B)\overline{\otimes}P]\rtimes_{\pi}G$ and let $Q$ be a diffuse von Neumann subalgebra of $(\overline{\bigotimes}_{L}B)\overline{\otimes}P$, for some subset $L$ of $ X$.
If no corner of $Q$ embedds into $P$(e.g. if $P=\Bbb C1$), then $$\Cal N_{M}(Q)\subset\sum_{g\in K}[(\overline{\bigotimes}_{X}B)\overline{\otimes}P]u_g,$$ where $K=\{g\in G|\exists x,y\in L, gx=y\}$ and $\{u_g\}_{g\in G}$ denote the canonical unitaries implementing  $\pi$.
\endproclaim
{\it Proof.} Let $w=\sum_g a_gu_g\in \Cal N_{M}(Q)$ and denote $\theta=Ad(w)_{\mid Q}\in \text{Aut}(Q)$, thus $$a_g\pi_g(y)=\theta(y)a_g,\forall g\in G,\forall y\in Q.$$

Now, let $a=(\otimes_{x\in X}a_x)\otimes u,b=(\otimes_{x\in X}b_x)\otimes v\in(\bigotimes_{X}B)\otimes P$, where $a_x,b_x\in B,$ for all $x\in X$, only finitely many different from $1$, and  $u,v\in\Cal U(P).$
We claim that for all $g\in G\setminus K$, there exists $a'\in M$, such that $$||E_{[(\overline{\bigotimes}_{L}B)\overline{\otimes}P]}(a\pi_g(y)b)||_2=||E_{P}(a'y)||_2,\forall y\in (\overline{\bigotimes}_{L}B)\overline{\otimes}P.$$
\vskip 0.05in
To see this, let $g\in G\setminus K$ (i.e. such that $gL\cap L=0$) and let $\{\xi_i|i\in I\}\subset L^2(P)$ be an orthonormal basis. 
Write $y\in (\overline{\bigotimes}_{L}B)\overline{\otimes}P$ as $y=\sum_{i\in I}y_i\otimes \xi_i$, for some $y_i\in L^2(\overline{\bigotimes}_{L}B)$. Then there exists  $a'',b''\in \overline{\bigotimes}_{gL}B$ and $c\in \overline{\bigotimes}_{X\setminus gL}B$ such that $$a\pi_g(y)b=\sum_{i\in I}  a''\sigma_g(y_i)b''\otimes c\otimes u\beta_g(\xi_i)v.$$ Thus,  $$||E_{[(\overline{\bigotimes}_{L}B){\overline{\otimes}}P]}(a\pi_g(y)b)||_2^2=
\sum_{i\in I}|\tau(c)|^2|\tau(a''\sigma_g(y_i)b'')|^2= ||E_{P}(a'y)||_2^2,$$
where $a'=|\tau(c)|$ $\sigma_{g^{-1}}({b''}^*a'')$.
\vskip 0.05in

Since no corner of $Q$ embedds into $P$, we can apply Popa's theorem (Theorem 1.1.)  to find $y_n\in\Cal U(Q),n\geq 1$  such that $$\lim_{n\rightarrow\infty}||E_{P}(a'y_n)||_2=0,\forall a'\in M.$$
Using the above claim and norm-$||.||_2$ approximations  we deduce that
$$\lim_{n\rightarrow\infty}||E_{[(\overline{\bigotimes}_{L}B)\overline{\otimes} P]}(a\sigma_g(y_n)b)||_2=0,\forall g\in G\setminus K,\forall a,b\in (\overline{\bigotimes}_{X}B)\overline{\otimes}P.$$

Combining this with  the  inequality $$  |\tau (a_g\sigma_g(y_n){{a_g}^*}\theta({{y_n}^*}))|\leq ||E_{[(\overline{\bigotimes}_{L}B)\overline{\otimes}P]}(a_g\sigma_g(y_n)a_g^*)||_2,$$ 

we get that $\lim_{n\rightarrow\infty} |\tau (a_g\sigma_g(y_n){{a_g}^*}\theta({{y_n}^*}))|=0,$ for all $ g\in G\setminus K.$
On the other hand, $\tau (a_g\sigma_g(y){{a_g}^*}\theta({{y}^*}))=\tau(a_g^*a_g),$ for all $y\in \Cal U(Q),$ thus
$a_g=0,$ for all $g\in G\setminus K$.\hfill $\blacksquare$
\vskip 0.05in

\proclaim{Corollary 1.6}
Let $N$ be a II$_1$ factor, $G$ be a countable discrete group and denote $M=N\wr G$. Suppose that $Q\subset M$ is a von Neumann subalgebra such that $\Cal N_{M}(Q)''$ is a factor and that a corner of $Q$ embedds into $\overline{\bigotimes}_{L}N$, for some finite set  $L\subset G$. Then there exists $u\in\Cal U(M)$ and $L',K\subset G$ finite such that  $uQu^*\subset \overline{\bigotimes}_{L}N$ and $u\Cal N_{M}(Q)''u^*\subset \sum_{g\in K}(\overline{\bigotimes}_{G}N)u_g$. Moreover, if $Q\subset \overline{\bigotimes}_{L}N$, then we can take $u=1$.
\endproclaim
{\it Proof.} If we denote $B=\overline{\bigotimes}_{L}N$, then both $B$ and $B'\cap M=\overline{\bigotimes}_{X\setminus L}N$ are factors. Thus, since $\Cal N_{M}(Q)''$ is a factor, by part $(i)$ of Corollary 1.3., we can find $u\in\Cal U(M)$ and $t>0$ such that $uQu^*\subset (\overline{\bigotimes}_{L}N)^t\subset (\overline{\bigotimes}_{L}N)^t\overline{\otimes}(\overline{\bigotimes}_{X\setminus L}N)^{1/t}$. 

If we let $g\in G\setminus L$ and denote $L'=L\cup\{g\}$, then since $N$ is a factor we can assume that $uQu^*\subset (\overline{\bigotimes}_{L}N)^t\overline{\otimes}N_{g}^{1/t}=\overline{\bigotimes}_{L'}N$. Finally, Lemma 1.5. implies that $u\Cal N_{M}(Q)''u^*\subset \sum_{g\in K}(\overline{\bigotimes}_{G}N)u_g$, where $K=L'{L'}^{-1}$.

\head 2. Weakly malleable actions. \endhead
The notion of {\it malleability} (see the definition below) has been introduced and used extensively by Popa in [Po2],[Po3],[Po5] and [Po6].
 Examples of malleable actions include both classical Bernoulli shift actions (i.e. of the form  $G\curvearrowright\overline{\bigotimes}_{G}L^{\infty}(X,\mu)$, where $(X,\mu)$ is a probability space) and non-commutative Bernoulli shift actions (i.e. of the form $G\curvearrowright\overline{\bigotimes}_{G}R$, where $R$ is  the hyperfinite II$_1$ factor or a matrix algebra, $\Bbb M_n(\Bbb C), 2\leq n\leq \infty$).

If we consider Bernoulli shift actions with the base $(N,\tau)$ an arbitrary finite von Neumann algebra, then  such actions are not malleable in general. However, we can  prove that Bernoulli shift actions  satisfy a  weaker version of Popa's malleability. In the next two sections we will show how results on malleable actions (e.g. "locating" rigid subalgebras of the associated crossed products, calculating 1-cohomology) can be extended to weakly malleable actions.

\vskip 0.05in
{\it Definition 2.1.} An action $\sigma$ of a countable discrete group $G$ on a finite von Neumann algebra $(N,\tau)$ is called  {\it weakly malleable } if there exists 
a finite von Neumann algebra $(\tilde N,\tilde\tau)$
 containing $N$ with $\tilde\tau_{|N}=\tau$, an action $\tilde\sigma: G\rightarrow$ Aut($\tilde N,\tilde\tau$) and a continuous action $\theta:\Bbb R\rightarrow$ Aut($\tilde N,\tilde\tau$) such that the following hold:
\vskip 0.05in
$(1)$. $\tilde\sigma(N)=N$ and $\tilde\sigma_{\mid N}=\sigma$.
\vskip 0.05in
$(2)$. $[\theta_t,\tilde\sigma(g)]=0,\forall t\in \Bbb R,\forall g\in G$.

\vskip 0.05in
$(3)$. There exists an orthonormal basis $\tilde\Cal N\subset \tilde N$ of $L^2(\tilde N)$ such that if $\tilde F\subset \tilde\Cal N$ is finite, then there exists $F\subset N$, $K\subset G$ finite sets  
 such that
 for all $x,y\in N, g\in {G\setminus K}$ and $\xi,\eta\in\tilde F$ we have that
$$|\tilde\tau(\theta_1(y^*)\xi x\sigma_g(\eta^*))|\leq
\sum_{\eta_1,\eta_2,\xi_1,\xi_2\in F}|\tau({\xi_1} y\sigma_g({\xi_2})||\tau({\eta_1}x\sigma_g({\eta_2}))|.$$
\vskip 0.05in
$(4)$. There exists $A\subset \tilde N$ a weakly dense $*$-algebra such that $$\lim_{g\rightarrow\infty}\sup_{x\in L^2(N),||x||_2\leq 1}||E_N((a-E_{N}(a))x\tilde{\sigma}_g(b-E_{N}(b)))||_2=0,\forall a,b\in A.$$ 
   
\vskip 0.1in

Moreover, if $\sigma$ is a mixing action, i.e.:
\vskip 0.05in
$(5)$. $\lim_{g\rightarrow\infty}\tau(x\sigma_g(y))=\tau(x)\tau(y),\forall x,y\in N,$
\vskip 0.05in
then is called {\it weakly malleable mixing}.
 \vskip 0.05in
 Note that condition (3) is equivalent to the existence of a finite set $F\subset N$ such that $$\tau(\theta_1(y^*)\xi x\sigma_g(\eta^*))=0,\forall  x,y\in N\ominus F\sigma_g(F), \forall g\in {G\setminus K},\forall \xi,\eta\in\tilde F.$$ Using this observation it follows that if $\sigma_i:G\rightarrow\text{Aut}(N_i,\tau_i)$ are two weakly malleable actions, then the product action $\sigma:G\rightarrow\text{Aut}(N_1\overline{\otimes}N_2,\tau_1\otimes \tau_2)$
given by $\sigma(g)={\sigma_1}(g)\otimes{\sigma_2}(g),$ for all $g\in G$, is also weakly malleable.
Also, remark that conditions (4) and (5) imply that $\tilde\sigma$ is mixing.

\vskip 0.1in
Recall from [Po2] that an action $\sigma$ of a discrete group $G$ on a finite von Neumann algebra $N$ is called  (tracial) {\it  malleable } if there exists 
 finite von Neumann algebra $\tilde N$
 containing $N$ an action $\tilde\sigma: G\rightarrow$ Aut($\tilde N,\tau$) and a continuous action $\theta:\Bbb R\rightarrow$ Aut($\tilde N,\tau$) such that conditions (1) and (2) hold true  as well as the following: 
\vskip 0.1in
$(6)$. $(\Cal N_{\tilde N}(N)\cap \theta_1(N))''=\theta_1(N).$
\vskip 0.1in
$(7)$. $\overline{\text{sp}}^{w}N\theta_1(N)=\tilde N$.
\vskip 0.1in
$(8)$. $\tau(x\theta_1(y))=\tau(x)\tau(y),\forall x,y\in N.$
\vskip 0.1in

\proclaim{Proposition 2.2} Any (tracial) malleable mixing action  is weakly malleable mixing.
\endproclaim
{\it Proof.} Let $\sigma:G\rightarrow \text{Aut}(N\subset \tilde N)$ be a tracial malleable mixing action and let $\theta:\Bbb R\rightarrow\text{Aut}(\tilde N)$  be the continuous action given by the definition of malleability. To prove that is weakly-malleable, we need to show that conditions (3) and (4) are satisfied.


Let $\{\eta_n\}_{n\geq 0}\subset N$ be an orthonormal basis for $L^2(N)$. Using condition (8), we get that $\Cal N=\{\theta_1(\eta_m)\eta_n\}_{m,n\geq 0}$ is an orthonormal family of $\tilde N$.	
Also, it is immediate that sp $\Cal N$ is $||.||_2$ dense in sp $\theta_1(N)N$. Since by condition (7) 
$\overline{\text{sp}}^{w}\theta_1(N)N=\tilde N$, we derive that sp $\Cal N$ is weakly dense in $\tilde N$.

This implies  that $\Cal N$ is an orthonormal basis for $L^2(\tilde N)$ and since $\forall m,n,p,q\geq 0,\forall x,y\in N$ and $\forall g\in G$ we have that: $$\tau(\theta_1(y^*)\theta_1(\eta_m)\eta_nx\sigma_g(\eta_q^*{\theta_1}(\eta_p^*)))=\tau(\sigma_g(\eta_p^*)y^*\eta_m)\tau(\eta_nx\sigma_g(\eta_q^*)),$$ it follows that
condition (3) holds true. 

To prove condition (4), let $\Cal G=\Cal N_{\tilde N}(N)\cap\theta_1(N)$. Then $A=$ sp $\Cal GN\subset\tilde N$ is a weakly dense $*$-subalgebra. 
Remark that if $u\in\Cal G$, then by projecting onto $N$ in the equality $\text{Ad}(u)(x)u=ux,\forall x\in N$, we get that $$\tau(u) [\text{Ad}(u)(x)-x)]=0,\forall x\in N.$$
Let $a=ux,b=vy\in A$ with $u,v\in\Cal G,x,y\in N$. Then for all $\xi\in L^2(N)$ and for all $g\in G$, $$E_N((a-E_{N}(a))\xi\sigma_g(b-E_{N}(b)))=E_{N}((ux-\tau(u)x)\xi\sigma_g(vy-\tau(v)y))=$$
$$\tau(u\sigma_g(v))[\text{Ad}(u)(x\xi)\sigma_g(y)]-\tau(u)\tau(v)[x\xi\sigma_g(y)]=$$
$$[\tau(u\sigma_g(v))-\tau(u)\tau(v)][\text{Ad}(u)(x\xi)\sigma_g(y)].$$

Since $\sigma$ is mixing we have that $\lim_{g\rightarrow\infty}|\tau(u\sigma_g(v))-\tau(u)\tau(v)|=0$. Combining this with the inequality $$||\text{Ad}(u)(x\xi)\sigma_g(y)||_2\leq ||x|| ||y|| ||\xi||_2,\forall\xi\in L^2(N),$$ we conclude that $\lim_{g\rightarrow\infty}\sup_{\xi\in L^2(N),||\xi||_2=1}||E_N((a-E_{N}(a))\xi\sigma_g(b-E_{N}(b)))||_2=0$ for all $a,b\in A$ of the above form. As such elements span $A$, the proof is complete.	
\hfill $\blacksquare$
\vskip 0.1in

\proclaim{Proposition 2.3}Let $G$ be a countable discrete group  and let $(B,\tau)$ be a finite von Neumann algebra.  Then
 the Bernoulli shift action $\sigma: G\rightarrow \text{Aut}(\overline{\bigotimes}_{G}B)$ is weakly malleable mixing.

\endproclaim
{\it Proof.} Define $\tilde B=B*L(\Bbb Z)$ and let $u\in L(\Bbb Z)$ be the canonical generating Haar unitary. Let $h=h^*$ be a hermitian element such that $u=\exp(ih)$ and set $u_t=\exp(ith),$ for all $t\in\Bbb R$.
  Let $N=\overline{\bigotimes}_{G}B$, $\tilde N=\overline{\bigotimes}_{G} \tilde B$ and let $\tilde\sigma:G\rightarrow\text{Aut}(\tilde N)$ be the Bernoulli shift action, which clearly extends $\sigma$.
Then  $\theta:\Bbb R\rightarrow\text{Aut}(\tilde N)$ given by $$\theta_t(\otimes_g{a_g})=\otimes_g\text{Ad}(u_t)(a_g),\forall t\in\Bbb R,\forall a={(a_g)}_g\in\tilde N$$ defines a continuous action of $\Bbb R$ on $\tilde N$ which  commutes with $\tilde\sigma$. Thus, to  prove that $\sigma$ is weakly malleable, we only need to check  conditions 3 and 4  in Definition 2.1, while mixingness is a well-known property of Bernoulli shifts.

Through the Gram-Schmidt procedure we can produce $\Cal B=\{1=\eta_0,\eta_1,..\}$ an orthonormal for $L^2(B,\tau)$ such that $\Cal B\subset B$ and that sp $\Cal B$ is a $*$-algebra (i.e. $\eta_i^*,\eta_i\eta_j\in$ sp $\Cal B,\forall i,j$).
Then an orthonormal basis for $\tilde B$ is given by $$\tilde{\Cal B}=\{u^{n_1}\eta_{j_1}u^{n_2}...\eta_{j_k}|j_1,..,j_{k-1},n_2,..,n_k\geq 1,k\in \Bbb N\}=\{\tilde{\eta}_n\}_{n\geq 0}.$$
Also, we have that
$\Cal N=\{\eta_I=\otimes_{g}{\eta_{i_g}}| |\{g\in 
G|i_g\not= 0\}|<\infty \},$ and
$\tilde{\Cal N}=\{\tilde{\eta}_I=\otimes_{g}{\tilde{\eta}_{i_g}}| |\{g\in G|i_g\not= 0\}|<\infty \}$
are orthonormal basis for $N$ and $\tilde N$, respectively. For $a=\otimes_g{a_g}\in \tilde N$, we define supp$(a)=\{g\in G|a_g\not= 1\}$ and if $S$ is a finite collection of such elements of $\tilde N$, then we set supp$(S)=\cup_{a\in S}$ supp$(a)$. 

To prove that condition 3 in definition 2.1. holds true, let $ \tilde F\subset \tilde{\Cal N}$ be finite.
Note that since $\Cal B^*\subset$ sp $B$, we have that $\theta_1(y^*)x,x^*\theta_1(y)\in$ sp $\tilde \Cal N,\forall x,y\in\Cal N$.
Thus,  there exists  $F\subset \Cal N$ finite such that if $x,y\in\Cal N,\eta\in \tilde F$, then $ \tau(\theta_1(y^*)x\eta^*)\not=0$ or  $\tau(\theta_1(y^*)\eta x)\not= 0$ implies $x,y\in F.$

  Let $K\subset G$ be a finite set such that for all $g\in G\setminus K$, $A=\text{supp}(\tilde F)$ and $B_g=\text{supp}(\sigma_g(\tilde F))$ are disjoint. Fix $g\in G\setminus K$ and set $C_g=G\setminus(A\cup B_g)$.
Given  $x,y\in\Cal N$, write $x=x_A\otimes x_{B_g}\otimes x_{C_g}$ and $ y=y_A\otimes y_{B_g}\otimes y_{C_g}$, where $x_D,y_D\in\Cal N$ are supported on $D$, for every set $D\in\{A,B_g,C_g\}$. Then, we have that $$\tau(\theta_1(y^*))\xi x\sigma_g(\eta^*))=\tau(\theta_1(y_{C_g}^*)x_{C_g})\tau(\theta_1(y_{B_g}^*)x_{B_g}\sigma_g(\eta^*))\tau(\theta_1(y_A^*)\xi x_A),\forall\xi,\eta\in\tilde F.$$ 
Thus, if  $\tau(\theta_1(y^*))\xi x\sigma_g(\eta^*))\not=0$,  for some $\xi,\eta\in\tilde F,x,y\in\Cal N$, then by using the freeness and the above we deduce that $$x_{C_g}=y_{C_g}=1,x_{B_g},y_{B_g}\in \sigma_g(F),x_A,y_A\in F,$$
hence that $x,y\in F\sigma_g(F)$.

Altogether, we have that $$\tau(\theta_1(y^*)\xi x\sigma_g(\eta^*))=0,$$
$$\forall g\in G\setminus K, \forall \xi,\eta\in \tilde F,\forall x,y\in  (\Cal N\times\Cal N)\setminus((F\sigma_g(F))\times (F\sigma_g(F)), $$
which proves (3).

Finally, note that (4) is verified for $A=$ sp $\tilde\Cal N$ (which is a $*$-algebra since sp $\Cal B$ is), since if $a,b\in\tilde \Cal N\setminus\Cal N$  then  the supports of $a$ and the  $\sigma_g(b)$ become disjoint, as $g\rightarrow \infty$, hence $aN\sigma_g(b)\perp N$, for $g$ outside a finite set.  \hfill $\blacksquare$

\vskip 0.05in
For the next lemma, assume the same context as in the proof of Proposition 2.3. and define $M=N\rtimes_{\sigma}G,\tilde M=\tilde N\rtimes_{\tilde\sigma}G$. Since $[\theta_t,\tilde\sigma_g]=0,$ for all $t\in\Bbb R$ and all $g\in G$, we have that $\theta_t$ extends to an automorphism of $\tilde M$ (still denoted $\theta_t$), which leaves invariant $M$.
 
\proclaim {2.4 Lemma}  Let $P\subset  M$ be a von Neumann subalgebra. Suppose that there exists $t>0$, a projection $p\in P,$ a non-zero  partial isometry $v\in \tilde M$   satisfying  $v^*v\leq p$ and $\theta_t(x)v=vx,$ for all $x\in pPp$. Then there exists $g\in G$, $F\subset G$ finite  and $c>0$ such that $$||E_{[\overline{\bigotimes}_{F}\tilde B]}(\theta_t(u)vu_g^*)||_2\geq c,\forall u\in\Cal U(pPp).$$
\endproclaim
{\it Proof.} Write $v=\sum_{g\in G}v_gu_g$, where $v_g\in\tilde N$ and let $g\in G$ be such that $v_g\not=0$. The hypothesis implies that $\theta_t(x)v_g=v_g\sigma_g(x),\forall x\in pPp.$

For $\varepsilon>0$, let $v_g'\in \tilde N$ such that $F=\text{supp}(v_g')$ is finite and $||v_g-v_g'||_2\leq \varepsilon$. 
Thus, $$||\theta_t(x)v_g'-v_g'\sigma_g(x)||_2\leq 2\varepsilon, \forall x\in (pPp)_1.$$
Note that $\theta_t(x)v_g'\in\Cal H=L^2([\overline{\bigotimes}_{G\setminus F}\theta_t(B)]\overline{\otimes}[\overline{\bigotimes}_{F}\tilde B]),$ for all $x\in pPp$. Hence, if we denote by $T$ the orthogonal projection onto $\Cal H$, then  $$||T(v_g'\sigma_g(x))-v_g'\sigma_g(x)||_2\leq 4\varepsilon,\forall x\in (pPp)_1.$$

Now, observe that if $x\in B$, then $E_{\theta_t(B)}(x)=|\tau(u_t)|^2\theta_t(x)$. Using this observation it is easy to see that if $\xi\in L^2([\overline{\bigotimes}_{G\setminus F}B]\overline{\otimes}[\overline{\bigotimes}_{F}\tilde B]),$
then $$||T(\xi)||_2^2\leq |\tau(u_t)|^4||\xi||_2^2+(1-|\tau(u_t)|^4)||E_{[\overline{\bigotimes}_{F}\tilde B]}(\xi)||_2^2.$$
Applying this inequality to $\xi=v_g'\sigma_g(x)$ with $x\in pPp$ and using the fact that $$||T(v_g'\sigma_g(x))||_2^2\geq ||v_g'\sigma_g(x)||_2^2-(4\varepsilon)^2, \forall x\in (pPp)_1,$$ we further get that $$(1-|\tau(u_t)|^4)||E_{[\overline{\bigotimes}_{F}\tilde B]}(v_g'\sigma_g(u))||_2^2\geq (1-|\tau(u_t)|^4)||v_g'\sigma_g(u)||_2^2-(4\varepsilon)^2\geq$$ $$(1-|\tau(u_t)|^4)(||v_g||_2-\varepsilon)^2-(4\varepsilon)^2,\forall u\in\Cal U(pPp).$$

If we choose $\varepsilon$ sufficiently small, then we get that there exists $c>0$ such that $$||E_{[\overline{\bigotimes}_{F}\tilde B]}(v_g\sigma_g(u))||_2\geq c,\forall u\in\Cal U(pPp).$$ Finally, since 
$$ E_{[\overline{\bigotimes}_{F}\tilde B]}(v_g\sigma_g(u))=E_{[\overline{\bigotimes}_{F}\tilde B]}(vuu_g^*)=E_{[\overline{\bigotimes}_{F}\tilde B]}(\theta_t(u)vu_g^*),\forall u\in \Cal U(pPp),$$ we are done.\hfill $\blacksquare$
\vskip 0.1in

\head 3. Relatively rigid subalgebras of wreath products \endhead
In this section we use the deformation properties of weakly malleable actions  to obtain conjugacy results for  rigid subalgebras of the associated crossed products. 
\vskip 0.1in
{\it Notation 3.1.} Throughout the next two sections, we fix a weakly malleable action  $\sigma: G\rightarrow \text{Aut}(N)$.
 Let $\tilde N\supset N$ and $\tilde\sigma:G\rightarrow \text{Aut}(\tilde N),$ $\theta:\Bbb R\rightarrow\text{Aut}(\tilde N)$ be as in Definition 2.1. Define $M=N\rtimes_{\sigma}G$, $\tilde M=\tilde N\rtimes_{\tilde\sigma}G$ and denote by $\{u_g\}_{g\in G}$ the canonical unitaries implementing the action of $G$ on $\tilde N$.  Since the action $\theta:\Bbb R\rightarrow\text{Aut}(\tilde N,\tau)$ commutes with $\tilde\sigma$, it extends to  a continuous action $\theta:\Bbb R\rightarrow\text{Aut}(\tilde M,\tau)$, given by $\theta_t(\sum_{g\in G}x_gu_g)=\sum_{g\in G}\theta_t(x_g)u_g.$

In the case $N=\overline{\bigotimes}_{G}B$ for some finite von Neumann algebra, $(B,\tau)$, and $\sigma:G\rightarrow\text{Aut}(N)$ is the Bernoulli $G$-action  we will consider the context from the proof of Proposition 2.3. More precisely, we let $\tilde N=\overline{\bigotimes}_{G} (B*L(\Bbb Z))$, $\tilde\sigma:G\rightarrow\text{Aut}(\tilde N)$ be the natural extension of $\sigma$ and $\theta_t\in\text{Aut}(\tilde N)$ be given by $\theta_t=\otimes_{g\in G}(\text{Ad} (u_t))_g,$ for all $t\in\Bbb R$, where $u_t=\exp(t\log(u))$, with $u\in L(\Bbb Z)$ being the generating Haar unitary.
\vskip 0.1in
We begin  this section by proving a technical result  about intertwiners in algebras associated to weakly malleable actions. This result is analogous to Theorem 3.2 in [Po2], which deals with the same problem in the case of malleable actions.
The proof is based on the following principle from section 3 in [Po2] : assume that $Q\subset N\subset  M$  are finite von Neumann algebras such that there exists $v_n\in\Cal U(Q)$ satisfying $\lim_{n\rightarrow\infty}||E_{N}(xv_ny)||_2=0,$ for all $x,y\in M\ominus N$. If $x\in M$ is so that $Qx\subset \sum_i x_iN$ for some $x_i\in  M$, then $x\in N$.
\vskip 0.1in

\proclaim{Proposition 3.2} Suppose that $Q\subset M$ is a von Neumann subalgebra such  that no corner of $Q$ embedds into $N$ inside $M$.
If $x\in\tilde M$ satisfies $Qx\subset \sum_i x_iM$ for some $x_i\in \tilde M$, then $x\in M$. In particular, $Q'\cap \tilde M\subset M$. 
 \endproclaim
{\it Proof.} Since no corner of $Q$ embbeds into $N$, we can find a sequence of unitaries $v_n\in\Cal U(Q),n\geq 1$ such that $$\lim_{n\rightarrow\infty}||E_{N}(v_nu_g^*)||_2=0,\forall g\in G.$$
We claim that  $$\lim_{n\rightarrow\infty}||E_{M}(xv_ny)||_2=0,\forall x,y\in \tilde M\ominus M,$$ 
which by the discussion preceding this proposition, gives us the conclusion.

Now, let $A\subset \tilde N$ be the weakly dense $*$-subalgebra for which condition (3) in Definition 2.1.  holds true. By Kaplansky's density theorem it suffices to verify the above claim for $x$ and $y$ of the form $x=(a-E_{N}(a))u_h,y=(b-E_{N}(b))u_k$, where $a,b\in A$  and  $h,k\in G$.
Set $a'=a-E_{N}(a),b'=b-E_{N}(b)$ and $$\alpha_g=\sup_{x\in L^2(N),||x||_2\leq 1} ||\sigma_{h^{-1}}(a')x\sigma_g(b')||_2.$$ Since $a,b\in A$, condition (3) in Definition 2.1. gives that $\lim_{g\rightarrow\infty}\alpha_g=0$.

 Next, observe that  $E_{M}(xv_ny)=\sum_{g\in G}E_{N}(a'\sigma_h(E_{N}(v_nu_g^*))\sigma_{hg}(b'))u_{hgk}$, thus $$||E_{M}(xv_ny)||_2^2=\sum_{g\in G}||\sigma_{h^{-1}}(a')E_{N}(v_nu_g^*)\sigma_g(b'	)||_2^2\leq$$
$$\sum_{g\in G}\alpha_g^2||E_{N}(v_nu_g^*)||_2^2\leq \max_{g\in L}\alpha_g^2\sum_{g\in L}||E_{N}(v_nu_g^*)||_2^2+\sup_{g\in G\setminus L}\alpha_g^2,$$
 for every finite set $L\subset G$. Altogether, we get that $\lim_{n\rightarrow\infty}||E_{M}(xv_ny)||_2=0$, for every $x,y$ as above, thus ending the proof.
\hfill $\blacksquare$
\vskip 0.05in

Recall that if $(N,\tau)$ is a finite von Neumann algebra, then a von Neumann subalgeba $Q\subset N$ is called  {\it relatively rigid} if any deformation of $id_N$ by subunital, subtracial, c.p. maps $\{\phi_n\}_n$ converges uniformly to $id_N$ on the unit ball of $Q$, i.e. $\lim_{n\rightarrow\infty}\sup_{x\in (Q)_1}||\phi_n(x)-x||_2=0$ ([Po1]). 

\proclaim{ Theorem 3.3}  Let $\sigma: G\rightarrow \text{Aut}(N)$ be a weakly malleable action and denote $M=N\rtimes_{\sigma}G$. Assume that $Q\subset M$ is a relatively rigid von Neumann subalgebra. Suppose that 
 
(i) $P=\Cal N_M(Q)''$ is a factor or that

(ii) $\sigma$ is a  Bernoulli shift action.

Then  a corner of $Q$ embbeds into $N$ or into $L(G)$.
\endproclaim

The proof of this theorem is splitted into the following two lemmas. Before proceding to the lemmas, we note that the proof of Theorem 3.3. can be easily adapted to show:
\vskip 0.05in
\proclaim {Theorem 3.3'}   Let $\sigma: G\rightarrow \text{Aut}(N)$ be a weakly malleable action and $\alpha: G\rightarrow\text{Aut}(N')$ be an action. Denote $M=(N\overline{\otimes}N')\rtimes_{\sigma\times\alpha}G$ and assume that $Q\subset M$ is a relatively rigid von Neumann subalgebra. Suppose that 
 
(i) $P=\Cal N_M(Q)''$ is a factor or that

(ii) $\sigma$ is a  Bernoulli shift action.

Then  a corner of $Q$ embbeds into $N\overline{\otimes}N'$ or into $N'\rtimes_{\alpha}G$.
\endproclaim

\proclaim {Lemma 3.4}  Assuming the context of 3.1. and 3.3., then one of the following is true:

$(a)$ There exists $w\in M,w\not=0$ such that $\theta_1(x)w=wx,$  for all $x\in Q$.

$(b)$ A corner of $Q$ embedds into $N$.
\endproclaim
{\it Proof.} Since $Q\subset M$ is rigid, $Q\subset\tilde M$ is also rigid([Po1]). Thus, we can find 
$F\subset \tilde M$ finite and $\delta>0$ such that if $\phi:\tilde M\rightarrow\tilde M$ is a normal, subunital, c.p. map with $||\phi(x)-x||_2\leq\delta,$ for all $x\in F$, then $||\phi(u)-u||_2\leq 1/2,$  for all $u\in\Cal U(Q)$.
In particular, since $t\rightarrow\theta_t$ is a pointwise $||.||_2$-continuous action, we can find $n\geq 1$ such that $||\theta_{1/2^n}(u)-u||\leq 1/2,\forall u\in\Cal U(Q).$ 

Let $v$ be the minimal $||.||_2$ element of $\Cal K=\overline{co}^{w}\{\theta_{1/2^n}(u)u^*|u\in\Cal U(Q)\}$. Then, since $||\theta_{1/2^n}(u)u^*-1||_2\leq 1/2,\forall u\in\Cal U(Q)$, we get that  $||v-1||_2\leq 1/2$, thus $v\not=0$. Also, since $\theta_{1/2^n}(u)\Cal Ku^*=\Cal K$ and $||\theta_{1/2^n}(u)xu^*||_2=||x||_2,$ for all $u\in\Cal U(Q),$ we deduce, using the uniqueness of $v$, that $\theta_{1/2^n}(u)v=vu,$ for all $u\in \Cal U(Q)$. 

Assume that $(b)$ is false, i.e. no corner of $Q$ embedds into $N$.
We then claim that since $\theta_{1/2^n}(u)v=vu,$ for all $u\in \Cal U(Q)$, then we can find $w\in\tilde M$ a non-zero partial isometry, such that $\theta_{1}(u)w=wu,$ for all $u\in \Cal U(Q)$, which proves $(a)$. To show this claim we treat separately the two cases.

$(i)$ $P$ {\it is a factor.}
Since no corner of $Q$ embedds into $N$, then, in particular, no corner of $P=\Cal N_{M}(Q)''
$ embedds into $N$. We can thus apply Proposition 3.1. to deduce that $P'\cap \tilde M\subset M$. As by our assumption, $P$ is a factor, or equivalently $P'\cap M=\Bbb C1$, we get that $P'\cap \tilde M=\Bbb C1$.
Now, by the proof of Theorem 4.1. in [Po2], steps 1-3, combined with $P'\cap \tilde M=\Bbb C1$ we get  the claim in the case $(i)$.

\vskip 0.05in
$(ii)$ $\sigma$ {\it is a Bernoulli action}. Let $\beta$ be the automorphism of $\tilde N$ given by $\beta_{|N}=id_{|N}$ and $\beta((u)_g)=(u)_g^*,$ for all $g\in G$. Then $\beta$ commutes with $\tilde \sigma$ and thus it extends to an automorphism of $\tilde M=\tilde N\rtimes_{\tilde \sigma}G$. Moreover, $\beta$ satisfies $\beta^2=id$, $\beta\theta_t\beta=\theta_{-t},$ for all $t\in \Bbb R$ and $M\subset {\tilde M}^{\beta}$. Finally, the same computations as in the proof of Proposition 3.3.in [IPP], which use only  the fact that $Q'\cap\tilde M\subset M$, quaranteed here by Proposition 3.2.(since no corner of $Q$ embedds into $N$), give us the claim in this case.
\hfill $\blacksquare$	
\vskip 0.1in
\proclaim{ Lemma 3.5}  Suppose that $Q\subset M=N\rtimes_{\sigma}G$ is a von Neumann subalgebra such that $\theta_1(x)w=wx,$ for all $x\in Q,$ for some non-zero partial isometry $w\in M$. Then a corner of $Q$ embedds into $N$ or into $L(G)$.

\endproclaim
{\it Proof.} Let $\tilde\Cal N\subset \tilde N$ be an orthonormal basis for which condition $(4)$ in Definition 2.1. holds. 
Then, let $e\in S=S^{-1} \subset G$ finite, $\tilde F\subset \tilde\Cal N$ finite  and $v=\sum_{(\xi,k)\in {\tilde F}\times S}\alpha_{\xi,k}\xi u_k$ such that $||v-w||_2\leq ||w||_2/3$, where $\alpha_{\xi,k}\in \Bbb C1$. Then it follows easily that $$ |\tau(\theta_1(u^*)vu{v}^*)|\geq ||w||_2^2/5>0,\forall u\in\Cal U(Q).$$ 

\vskip 0.05in
{\it Claim.} There exists $a_i\in N,i=\overline{1,n}$ such that $$|\tau(\theta_1(u^*)vu{v}^*)|\leq \sum_{i,j}(||E_{N}(a_iu{a_j}^*)||_2^2+||E_{L(G)}(a_iu{a_j}^*)||_2^2),\forall u\in\Cal U(N).$$
\vskip 0.05in
Note that this claim combined with the above inequality gives that $$||w||_2^2/5\leq \sum_{i,j}(||E_{N}(a_iu{a_j}^*)||_2^2+||E_{L(G)}(a_iu{a_j}^*)||_2^2),\forall u\in\Cal U(N).$$ However, if no corner of $Q$ could be embedded neither into $L(G)$ nor into $N$, then by Popa's theorem (Theorem 1.1., see also [IPP]), given any $a_1,..,a_n\in M,\varepsilon>0$, we can find $u\in\Cal U(Q)$ such that $||E_{L(G)}(a_iua_j^*)||_2,||E_{N}(a_iua_j^*)||_2\leq\varepsilon,\forall i,j$, a contradiction.
\vskip 0.05in
To prove the claim, let $u\in\Cal U(N)$ and write $u=\sum_{g\in G}x_gu_g$, where $x_g=E_N(x{u_g}^*)$. Then a direct computation shows that $$\tau(\theta_1(u^*)vu{v}^*)=\sum_{g\in G, (\xi,k),(\eta,l)\in {\tilde F}\times S}
\alpha_{\xi,k}\overline{\alpha_{\eta,l}}\tau(\theta_1({x_g}^*)\xi\sigma_k(x_{k^{-1}gl})\sigma_{g}({\eta}^*)).$$
By Definition 2.1.(4), we can find $F\subset N,K\subset G$ finite sets such that for all $x,y\in N, g\in {G\setminus K}$ and $ \xi,\eta\in\tilde F$ we have that
$$|\tau(\theta_1(y^*)\xi x\sigma_g(\eta^*))|\leq
\sum_{\eta_1,\eta_2,\xi_1,\xi_2\in F}|\tau(\xi_1 y\sigma_g(\xi_2)||\tau(\eta_1x\sigma_g(\eta_2))|.$$
Thus, if we let $C=\max_{(\xi,k)\in {\tilde F}\times S}
|\alpha_{\xi,k}|^2$, then 
$$|\tau(\theta_1(u^*)vu{v}^*)|\leq C\sum_{g\in G\setminus K,k,l\in S}\sum_{\eta_1,\eta_2,\xi_1,\xi_2\in F}|\tau(\xi_1 x_g\sigma_g(\xi_2)||\tau(\eta_1\sigma_k(x_{k^{-1}gl})\sigma_g(\eta_2))|+$$
$$C\sum_{g\in K, (\xi,k),(\eta,l)\in {\tilde F}\times S}|\tau(\theta_1({x_g}^*)\xi\sigma_k(x_{k^{-1}gl})\sigma_{g}({\eta}^*))|.$$

Next, we estimate the first term:
$$\sum_{g\in G\setminus K,k,l\in S}\sum_{\eta_1,\eta_2,\xi_1,\xi_2\in F}|\tau(\xi_1 x_g\sigma_g(\xi_2)||\tau(\eta_1\sigma_k(x_{k^{-1}gl})\sigma_g(\eta_2))|\leq$$ 
$$1/2\sum_{g\in G,k,l\in S}\sum_{\eta_1,\eta_2,\xi_1,\xi_2\in F} (|\tau(\xi_1 x_g\sigma_g(\xi_2)|^2+|\tau(\eta_1 \sigma_k(x_{k^{-1}gl})\sigma_g(\eta_2))|^2)=$$
$$1/2\sum_{k,l\in S}(||E_{L(G)}(\xi_1u\xi_2)||_2^2+||E_{L(G)}(\sigma_{k^{-1}} (\eta_1)u\sigma_{{l^{-1}}(\eta_2}))||_2^2.$$
	
Similarly, for the second term we have:
$$\sum_{g\in K, (\xi,k),(\eta,l)\in {\tilde F}\times S}|\tau(\theta_1({x_g}^*)\xi\sigma_k(x_{k^{-1}gl})\sigma_{g}({\eta}^*))|$$
$$\leq 1/2\max_{\xi\in\tilde F}||\xi||^2 \sum_{g\in K,k,l\in S}(||x_g||_2^2+||x_{k^{-1}gl}||_2^2)=$$
 $$1/2\max_{\xi\in\tilde F}||\xi||^2 \sum_{g\in K,k,l\in S}(||E_N(u{u_g}^*)||_2^2+||E_N(uu_{k^{-1}gl^*})||_2^2).$$
Altogether, we get: $$|\tau(\theta_1(u^*)vu{v}^*)|\leq C \sum_{k,l\in S,\xi_1,\xi_2\in F}
||E_{L(G)}(\sigma_{k}({\xi_1}^*)u\sigma_{l}({\xi_2}^*))||_2^2+$$
$$C|SKS|\max_{\xi\in\tilde F}||\xi||^2\sum_{g\in SKS}||E_N(u{u_g}^*)||_2^2,\forall u\in \Cal U(N),$$
thus proving the claim.

\proclaim {Theorem 3.6 } Let $(B,\tau)$ be a finite von Neumann algebra  and let $G$ be a countable discrete group. Denote $M=B\wr G$ and suppose that $Q\subset M$ is a relatively rigid diffuse von Neumann subalgebra. Then either
\vskip 0.05in
$(i)$ A corner of $Q$ embedds into $L(G)$. Moreover, in this case, if $P=\Cal N_{M}(Q)''$ is a factor, then we can find $u\in\Cal U(M)$ such that $uQu^*\subset uPu^*\subset L(G)$.
\vskip 0.05in
$(ii)$ A corner of $Q$ embedds  into ${\overline{\bigotimes}_{L}B},$ for some finite set $L$. Moreover, in this case, if 
$P=\Cal N_{M}(Q)''$ and $B$ are factors, then we can find $L',K\subset G$ finite  and a unitary $u\in M$ such that $uQu^*\subset {\overline{\bigotimes}_{L'}B}$ and $uPu^*\subset \sum_{g\in K}(\overline{\bigotimes}_{G}B)u_g$.
   
\endproclaim
{\it Proof.} Following Theorem 3.3. we are in one of the following two cases:
\vskip 0.05in
 {\it (i) A corner of $Q$ embedds into $L(G)$.} 
\vskip 0.05in
If $P$ is a factor, then since  Bernoulli shift actions are properly outer and mixing, Theorem 3.1. in [Po2] implies that there exists $u\in\Cal U(M)$ such that $uPu^*\subset L(G)$.
\vskip 0.05in
 {\it (ii) A corner of $Q$ embedds into $\overline{\bigotimes}_{G}B$.}
\vskip 0.05in
In this case we only need to show that the first assertion, i.e. that a corner of $Q$ embedds into $\overline{\bigotimes}_{L}B$, for some finite set $L\subset G$. The second assertion will then follow directly from Corollary 1.6.

Since a corner of $Q$ embedds into $\overline{\bigotimes}_{G}B$, we can find a von Neumann subalgebra $P\subset\overline{\bigotimes}_{G}B$, projections $q\in Q,p\in P$, an  isomorphism $\phi:pPp\rightarrow qQq$ and a non-zero partial isometry $v$ such that $vv^*\leq q$ and $\phi(x)v=vx,$ for all $x\in pPp$. 

On the other hand, using the relative rigidity of $Q$ as in the proof of Lemma 3.4. we deduce that
if $w_t$ is the minimal norm-$||.||_2$ element of $\Cal K_t=\overline{co}^{w}\{\theta_{t}(u)u^*|u\in\Cal U(qQq)\}$, then $\lim_{t\rightarrow 0}w_t=q$ and $\theta_t(y)w_t=w_ty,$ for all  $y\in qQq$ and $t>0$. 
 Since $\theta_t\rightarrow$ id as $t\rightarrow 0$, we have that $\lim_{t\rightarrow 0}\theta_t(v^*)w_tv=v^*qv=v^*v$. In particular, we can find $t>0$ such that $\theta_t(v^*)w_tv\not=0$. Moreover, we have  that $$\theta_t(x)[\theta_t(v^*)w_tv]=[\theta_t(v^*)w_tv]x,\forall x\in pPp.$$
By Lemma 2.4., which applies since $P\subset \overline{\bigotimes}_{G}B$, we can find $g\in G$, $F\subset G$ finite and $c>0$ such that $$||E_{[\overline{\bigotimes}_{F}\tilde B]}(\theta_t(uv^*)w_tvu_g^*)||_2\geq c,\forall u\in\Cal U(pPp).$$
Since $uv^*=v^*\phi(u),$ for all $u\in\Cal U(pPp)$, $\phi(\Cal U(pPp))=\Cal U(qQq)$ and $\overline{\bigotimes}_{F}\tilde B$ is invaried by $\theta_t$, the last inequality is equivalent to $$||E_{[\overline{\bigotimes}_{F}\tilde B]}(v^*U \theta_{-t}(w_tv)u_g^*)||_2\geq c,\forall U\in\Cal U(qQq).$$
Using the fact that $Q\subset M$, this is  easily  seen to imply that we can find $a_1,a_2,..,a_n\in M$ and $L\subset G$ finite such that $$\sum_{i,j}||E_{[\overline{\bigotimes}_{L}B]}(a_iUa_j^*)||_2\geq c/2,\forall U\in\Cal U(qQq),$$ which by Popa's theorem implies that a corner of $Q$ embedds into $\overline{\bigotimes}_{L}B$.
\hfill $\blacksquare$
\vskip 0.1in
\proclaim{Corollary 3.7} Let $(B,\tau)$ be a finite von Neumann algebra having Haagerup's property and let $G$ be a countable discrete group. Denote $M=B\wr G$ and assume that $Q\subset M$ is a relatively rigid diffuse von Neumann subalgebra such that $P=\Cal N_{M}(Q)''$ is a factor. Then there exists $u\in\Cal U(M)$ such that $uQu^*\subset uPu^*\subset L(G)$.
 
\endproclaim
{\it Proof.} By the last theorem we just have to show that no corner of $Q$ can be embedded into $\overline{\bigotimes}_{L}B$ for some finite set $L\subset G$. If we assume the contrary, then by Popa's theorem we can find $q\in\Cal P(Q),p\in\Cal P(\overline{\bigotimes}_{L}B)$, an isomorphism $\theta:qQq\rightarrow p(\overline{\bigotimes}_{L}B)p$ and a non-zero partial isometry $v\in M$ such that $\theta(x)v=vx,$ for all$ x\in qQq$ and $v^*v\leq q$.

Since $B$ has Haagerup's property, we can find unital, tracial compact c.p. maps $\phi_n:B\rightarrow B$ such that $\phi_n$ converges to {\it id}$_B$ pointwise. Then $\Phi_n:M\rightarrow M$ defined by $\Phi_n=\otimes_{g\in G}(\phi_n)_g$ on $\overline{\bigotimes}_{G}B$ and by $\Phi_n(u_g)=u_g$ for  $g\in G$, gives a deformation of ${id}_{M}$ by unital, tracial, c.p. maps. Denote by $T_{\Phi_n}:L^2M\rightarrow L^2M$ the induced bounded operator. Then $L^2(\overline{\bigotimes}_{L}B)$ is an invariant subspace for $T_{\Phi_n}$ and  $\Phi_{n_{|L^2(\overline{\bigotimes}_{L}B)}}$ is a compact operator.

Since $Q\subset M$ is rigid and $\Phi_n\rightarrow {id}_{M}$ pointwise, we can find $N\in\Bbb N$ such that $$||\Phi_n(x)-x||_2\leq \sqrt{\tau(q)}/3,$$ for all  $x\in (qQq)_1$ and all $n\geq N$.
Also, since $\Phi_n\rightarrow {id}_{M}$,  by Corollary 1.1.2 in [Po1], we get that $$\lim_{n\rightarrow\infty}\sup_{x\in (qQq)_1}||\Phi_n(\theta(x)v)-\Phi_n(\theta(x))v||_2=\lim_{n\rightarrow\infty}\sup_{x\in (qQq)_1}||\Phi_n(vx)-v\Phi_n(x)||_2=0.$$ Thus, using the equality $\theta(x)v=vx,$ for all $ x\in qQq$, we can find $n\geq N$ such that for all $x\in (qQq)_1$ we have that $$||\Phi_n(\theta(x))v-v\Phi_n(x)||_2\leq\sqrt{\tau(q)}/3.$$Combining this with the inequality
 $$||v\Phi_n(u)||_2\geq ||vu||_2-||v\Phi_n(u)-vu||_2\geq ||v||_2-||\Phi_n(u)-u||_2\geq$$
$$\sqrt{\tau(q)}-\sqrt{\tau(q)}/3=2\sqrt{\tau(q)}/3,\forall u\in\Cal U(qQq),$$
we obtain that $$||\Phi_n(\theta(u))v||_2\geq\sqrt{\tau(q)}/3,\forall u\in\Cal U(qQq).$$
 
Now, using the fact that $Q$ is diffuse we can find $u_m\in\Cal U(qQq),m\geq 1$ such that $u_m$ converges weakly to 0, thus $\theta(u_m)$ converges weakly to 0. Since $\theta(u_m)\in\overline{\bigotimes}_{L}B$ and ${T_{\Phi_n}}_{|L^2(\overline{\bigotimes}_{L}B)}$ is compact, we deduce that $\lim_{m\rightarrow\infty}||\Phi_n(\theta (u_m))||_2=0$, which contradicts the last inequality.
\vskip 0.1in
We end  this section by mentioning that a more general version of the above result holds true,  with its proof following closely the arguments used above. 
\proclaim{ Corollary 3.7'} Let $(B,\tau)$ be a finite von Neumann algebra having Haagerup's property and let $G$ be a countable discrete group. Let $\beta:G\rightarrow\text{Aut}(C)$ be an action of $G$ on a finite von Neumann and define $\alpha:G\rightarrow\text{Aut}((\overline{\bigotimes}_{G}B)\overline{\otimes} C)$ to be the product between the $(B,\tau)$-Bernoulli action and $\beta$. Denote $M=[(\overline{\bigotimes}_{G}B)\overline{\otimes} C]\rtimes_{\alpha}G$. If $Q\subset M$ is a relatively rigid von Neumann subalgebra such that $\Cal N_{M}(Q)''$ is a factor, then a corner of $Q$ embedds into $C\rtimes_{\beta}G$.
\endproclaim

\vskip 0.1in

\vskip 0.1in
\head 4. Applications\endhead
\proclaim{Theorem 4.1} 
Let $G,H$ be two discrete countable ICC   groups, let $P$ and $Q$ be two {\it weakly rigid} II$_1$ factors and denote $M=P\wr G$, $N=Q\wr H$.  
If $\theta:M\rightarrow N^t$ is an isomorphism for some $t>0$, then there exists $u\in\Cal U(N)$ such that $u\theta(\overline{\bigotimes}_{G}P)u=(\overline{\bigotimes}_{H}Q)^t$.
\endproclaim
{\it Proof.} Assume that $P\wr G=(Q\wr H)^t$ and denote by $u_g,v_h$ the unitaries implementing the actions of $G$ and $H$, respectively. If $Q_0\subset Q$ be a relatively rigid diffuse subalgebra such that $\Cal N_{Q}(Q_0)''=Q$, then it is clear that $ (\overline{\bigotimes}_{H}Q)^t\subset\Cal N_{(Q\wr H)^t}(Q_0^t)''$. Now, since $Q$ is a factor, we get that $(\overline{\bigotimes}_{H}Q)'\cap (Q\wr H)=\Bbb C1$, hence $[(\overline{\bigotimes}_{H}Q)^t]'\cap (Q\wr H)^t=\Bbb C1$. Combining the last two conclusions we deduce that $\Cal N_{(Q\wr H)^t}(Q_0^t)''$ is a factor.
Thus, we can apply Theorem 3.6. to the rigid inclusion $Q_0^t\subset M=P\wr G$, to get $ u\in \Cal U(M)$ and $K\subset G$ finite such that either $(i)$ $ u\Cal N_{(Q\wr H)^t}(Q_0^t)''u^*\subset L(G)$ or $ (ii)$ $ u\Cal N_{(Q\wr H)^t}(Q_0^t)''u^*\subset \sum_{g\in K}(\overline{\bigotimes}_{G}P)u_g$.

 Next, we argue that $(i)$ leads to a contradiction; if we assume $(i)$  then $\tilde Q:=u(\overline{\bigotimes}_{G}Q)^tu^*\subset L(G)$. Note that $\tilde Q$ is regular in $M$. On the other hand, since $\tilde Q\subset L(G)$,  Theorem 3.1. in [Po2]  gives that $\Cal N_{M}(\tilde Q)\subset L(G)$, a contradiction.

Therefore, it follows that $(ii)$ holds true, hence $u(\overline{\bigotimes}_{H}Q)^t u^*\subset \sum_{g\in K}(\overline{\bigotimes}_{G}P)u_g$ for $K\subset G$ finite. Similarly, we can find $v\in\Cal U(M)$ and $L\subset H$ finite such that $v(\overline{\bigotimes}_{G}P)v^*\subset\sum_{h\in L}(\overline{\bigotimes}_{H}Q)^t v_h$.
Thus, we get  two $\overline{\bigotimes}_{G}P-(\overline{\bigotimes}_{H}Q)^t$ Hilbert bimodules $\Cal H,\Cal K\subset L^2(M)$, which are finite dimensional over $\overline{\bigotimes}_{G}P$ and over $(\overline{\bigotimes}_{H}Q)^t$, respectively. 

Finally, since $G$ and $H$ are ICC groups, Theorem 8.4. in [IPP] can be  used to derive the conclusion.
\hfill $\blacksquare$
\vskip 0.1in
In [Po5], using deformation/rigidity techniques, the 1-cohomology of Connes-Stormer Bernoulli shifts is calculated. In particular, it is shown that  1-cocycles for the non-commutative Bernoulli shift actions of a weakly rigid group are cohomologous to characters. Inspired by the results and ideas in [Po5] we show that the same result holds for Bernoulli action with arbitrary base.

\vskip 0.05in
Let $\sigma$ be an action of a group $G$ on a II$_1$ factor $N$. Let $p\in\Cal P(N)$ be a non-zero projection and let $w:G\rightarrow N$ be a map such that $w_e=p$ and $w_g$ is a partial isometry with $w_gw_g^*=p,w_g^*w_g=\sigma_g(p),$ for all $g\in G$. If $w_{gh}=w_g\sigma_g(w_h),$ for all $g,h\in G$ then $w$ is called a {\it generalized 1-cocycle for} $\sigma$. 
\proclaim{Proposition 4.2} Let $\sigma:G\rightarrow\text{Aut}(N)$ be a weakly malleable mixing action and let $w:G\rightarrow N$ be a generalized 1-cocycle for $\sigma$. Suppose that $G_0\subset G$ is an infinite subgroup such that there exists a non-zero partial isometry $v\in \tilde N$ such that $\theta_1(p)vp=v$ and that $\theta_1(w_g)\sigma_g(v)=vw_g,$ for all $g\in G_0$.
 Then $p=1$ and  there exists $u\in\Cal U(N)$, $\gamma:G_0\rightarrow\Bbb T$ a character such that $w_g=\gamma_g u\sigma_g(u^*),$ for all $g\in G_0$. Moreover, if $G_0$ is normal in $G$, then $w_g=\gamma_g u\sigma_g(u^*),$ for all $g\in G,$
 for some unitary $u\in N$ and a character $\gamma:G\rightarrow\Bbb T$.
\endproclaim
{\it Proof}. 
 We first claim that  $v\in\Cal H=\overline{\text{sp}}$ $\theta_1(N)N$. 

{\it Proof of claim.}
 Since $\Cal H$ is a $\theta_1(N)-N$ bimodule, the relation $\theta_1(w_g)\sigma_g(v)=vw_g,$ for all $g\in G_0$ implies that $$\theta_1(w_g)\sigma_g(v-P_{\Cal H}(v))=(v-P_{\Cal H}(v))w_g,\forall g\in G_0.$$ Denote $v'=v-P_{\Cal H}(v)$, then $\theta_1(p)v'p=v'$ and $\theta_1(w_g)\sigma_g(v')=v'w_g,$ for all $g\in G_0$. Using these identities we get that $$\tau(\theta_1({w_g}^*)v'w_g\sigma_g({v'}^*))=\tau({v'}^*v'),\forall g\in G_0.$$
On the other hand, since $\sigma$ is weakly malleable,  we can find $F\subset N,K\subset G$ finite sets and $c>0$ such that for all $x,y\in N$, $g\in G\setminus K$ we have that:
	
$$|\tau(\theta_1(y^*)v' x\sigma_g({v'}^*))|\leq ||v'||_2^2/2+ c\max_{x,y\in F\sigma_g(F)}|\tau(\theta_1(y^*)v' x\sigma_g({v'}^*))|.$$
Thus, it follows that $$\max_{x,y\in F\sigma_g(F)}|\tau(\theta_1(y^*)v' x\sigma_g({v'}^*))|\geq ||v'||_2^2/(2c),\forall g\in G_0\setminus K.$$
Now, if $x_1,x_2,y_1,y_2\in F$ and $g\in G$, then $$\tau(\theta_1((y_1\sigma_g(y_2))^*)v'(x_1\sigma_g(x_2))\sigma_g({v'}^*))=\tau((\theta_1({y_1}^*)v'x_1) \sigma_g(x_1{v'}^*\theta_1({y_2}^*))).$$ 
Since $v'\perp \Cal H$, we get in particular that $\tau(\theta_1(y^*)v'x)=0,$ for all $x,y\in F$.
Combining this with the fact that $\sigma$ is mixing we derive that $$\lim_{g\rightarrow\infty} \tau((\theta_1({y_1}^*)v'x_1) \sigma_g(x_1{v'}^*\theta_1({y_2}^*)))=0.$$ This further implies that $\lim_{g\rightarrow\infty}\max_{x,y\in F\sigma_g(F)}|\tau(\theta_1(y^*)v' x\sigma_g({v'}^*))|=0$, hence $v'=0$, proving that $v\in\Cal H$.
\vskip 0.1in
In the second part of the proof we use an argument due to Popa to show that $p=1$ and that $w_{|G_0}$ is cohomlogous to a character. 
Recall from Proposition 3.4. in [Po5], that $$\sigma^w_g(\xi)=w_g\sigma_g(\xi),\forall \xi\in L^2(pN)$$ gives a unitary representation of $G$ on $L^2(pN)$. Also, if  we denote by $\Cal H\Cal S$ the Hilbert space of Hilbert-Schmidt operators on $L^2(pN)$, then $$\tilde{\sigma}_g^w(X)=\sigma_g^wX{\sigma_g^w}^*,\forall X\in\Cal H\Cal S$$ gives a unitary representation of $G$ on $\Cal H\Cal S$.

Note that the claim implies that $v$  lies in the Hilbert space $\Cal K=\overline{\text{sp}}$ $[\theta_1(pN)(Np)]$ and  satisfies $\theta_1(w_g)vw_g^*=v,$ for all $g\in G_0$.
Next we identify $\Cal K$ with  $\Cal H\Cal S$ in the natural way, by letting  $T(\theta_1(x)y)=L_{x,y}$, where $L_{x,y}\in\Cal H\Cal S$ is given by $L_{x,y}(\xi)=<\xi,y^*>x,$ for all $x\in pN, y\in Np$ and $ \xi\in L^2(pN)$. Since $$<\theta_1(x)y,\theta_1(x')y'>=\tau({x'}^*x)\tau(y{y'}^*)=Tr(L_{x',y'}^*L_{x,y}),\forall x,x'\in pN,\forall y,y'\in Np,$$ it is immediate  that $T$ extends to a unitary $T:(\Cal K,< >_{\tau})\rightarrow(\Cal H\Cal S,< >_{Tr})$.

Moreover, if $\xi\in \Cal K$ and $g\in G$, then $T(\theta_1(w_g)\sigma_g(\xi)w_g^*)=\tilde{\sigma}_g^w(T(\xi)).$
Indeed, if $x\in pN,y\in Np$, then it is easy to check that $$
\tilde{\sigma}_g^w(T(\theta_1(x)y)))=\tilde{\sigma}_g^w(L_{x,y})=L_{w_g\sigma_g(x),\sigma_g(y)w_g^*}=$$ $$T(\theta_1(w_g)\sigma_g(\theta_1(x)y)w_g^*).$$
Using the last equality, we deduce that $T(v)\in\Cal H\Cal S$ is a fixed vector for the representation $\tilde{\sigma}^w_{|G_0}$.
Proposition 3.4. in [Po5] then implies that there exists a finite dimensional Hilbert space $\Cal H_0\subset pN$ which is invariant to $\sigma^w_{|G_0}$. Moreover, the finite dimensional vector space $\Cal H_1=$ sp $\Cal H_0^*\Cal H_0$ is invariant under $\sigma_{|G_0}$.

Since $\sigma$ is a mixing action and $G_0$ is infinite, $\sigma_{|G_0}$ is weakly mixing, therefore we must have that $\Cal H_1=\Bbb C1$.
In particular, if $\xi\in \Cal H_0$ satisfies $||\xi||_2=1$, then $\xi^*\xi=1$, i.e. $\xi$ is a unitary. Thus, $\xi$ is a unitary whose left support lies under $p$, hence $p=1$. Also, since $\xi$ is a unitary and $\xi^*\Cal H_0\subset\Cal H_1=\Bbb C1$, we must have that $\Cal H_0=\Bbb C\xi$. This implies that there exists a character $\gamma:G_0\rightarrow\Bbb T$ such that $\sigma^w_g(\xi)=w_g\sigma_g(\xi)=\gamma(g)\xi,$ for all $g\in G_0$.

Thus, $w_g=\gamma(g)\xi\sigma_g(\xi^*),$ for all $g\in G_0$. If $G_0$ is normal in $G$, then using the mixingness of $\sigma$ it is immediate that $\xi^*w_g\sigma_g(\xi)\in \Bbb C1,$ for all $	g\in G$, which proves the conclusion.
\hfill $\blacksquare$
\vskip 0.05in
\proclaim {Corollary 4.3} Let $(B,\tau)$ be a finite von Neumann algebra, $G$ a countable discrete group and $\sigma:G\rightarrow\text{Aut}(\overline{\bigotimes}_{G}B)$ be the associated Bernoulli shift action. Let $w:G\rightarrow\overline{\bigotimes}_{G}B$ be a generalized 1-cocycle. 

$(i)$ If $G$ contains an infinite relatively rigid subgroup, then $p=1$.

$(ii)$ If $G$ contains an infinite normal relatively rigid subgroup, then H$^1(\sigma;G)=\text{Char}(G)$.
\endproclaim

{\it Proof.} Denote $M=B\wr G$.
 By the hypothesis there exists $G_0\subset G$ infinite subgroup such that the pair $(G,G_0)$ has relative property (T). Then $L(G_0)\subset L(G)$ is a rigid inclusion ([Po1]), thus $Q=\{w_gu_g|g\in G_0\}''\subset \{w_gu_g|g\in G\}''$ is a rigid inclusion.

This implies that $Q\subset pMp$ is a rigid subalgebra and by applying Lemma 3.4. (although $Q$ is not unital we can still apply Lemma 3.4. by working with suitable  amplifications) we get that either there exists $v\in M,v\not=0$ such that $\theta_1(p)vp=v$ and $\theta_1(x)v=vx,$ for all $x\in Q$ or a corner of $Q$ embedds into $\overline{\bigotimes}_{G}B$.
Since $\lim_{g\rightarrow\infty}||E_{\overline{\bigotimes}_{G}B}(aw_gu_gb)||_2=0,$ for all $a,b\in \tilde M$, we get that condition $(iii)$ in Theorem 1.1. holds true, thus no corner of $Q$ embedds into $\overline{\bigotimes}_{G}B$.

Thus $\theta_1(x)v=vx,$ for all $x\in Q$, which  for  $x=w_gu_g,g\in G_0$, gives that $\theta_1(w_g)\sigma_g(p)=pw_g,$  for all $g\in G_0$. We can now apply Proposition 4.2. to finish the proof. 

\proclaim{Corollary 4.4}
Let $G$ be a group having an infinite relatively rigid subgroup and $N$ be a weakly rigid II$_1$ factor. Then $\Cal F(N\wr G)=\{1\}$.
\endproclaim
{\it Proof.}  Denote $\tilde N=\overline{\bigotimes}_{G}N$ and $M=N\wr G$. Let $p\in\Cal P(\tilde N)$ and let $\alpha: M\rightarrow pMp\subset M$ be an isomorphism.  Using
Theorem 4.1. we can assume that $\alpha(\tilde N)=p\tilde Np $. This further implies that we can find partial isometries $w_g\in M$ such that $w_gw_g^*=p,w_g^*w_g=\sigma_g(p)$ and that $\alpha(u_g)=w_gu_g,$ for all $g\in  G.$  
Note that the relation $\alpha(u_{gh})=\alpha(u_g)\alpha(u_h),$ for all $g,h\in G$ implies that $w_g$ is a generalized 1-cocycle.

Since $G$ contains an infinite relatively rigid subgroup, by Corollary 4.3. we have that the support $p$ of $w$ must equal $1$, thus $\Cal F(N\wr G)=\{1\}$.
\hfill $\blacksquare$
\vskip 0.1in
{\it 4.5. Free products.} Let $(B,\tau)$ be a finite von Neumann algebra and $G$ be a countable discrete group. Define $N=*_{g\in G}B$ to be the free product of infinitely many copies of $B$ indexed by $G$. Then, similarly to Proposition 2.3., it follows that the {\it free Bernoulli} action $\sigma:G\rightarrow\text{Aut}(N)$ given by $\sigma_g(*_{h}x_h)=*_{h}x_{g^{-1}h},$ for all $g\in G$ is weakly malleable mixing. Using an argument paralleling the proof of Theorem 3.6., we can prove that if $Q\subset N\rtimes_{\sigma}G$ is a relatively rigid subalgebra whose normalizer generates a factor, then we can unitarily conjugate $Q$  into $L(G)$ or into $*_{F}B$, for some finite set $F\subset G$.
\vskip 0.05in

 Let $G_1,G_2$ be two countable discrete groups. We claim that if $Q\subset M=L(G_1*G_2)$ is a relatively rigid diffuse subalgebra such that $\Cal N_{M}(Q)''$ is a factor, then there exists $u\in\Cal U(M)$ and $i\in\{1,2\}$ such that $uQu^*\subset L(G_i)$ (equivalently, by [IPP], a corner of $Q$ embedds into $L(G_i)$). Note that this reproves a particular case of Theorem 4.3. in [IPP].

Let $Q\subset M$ be as above and assume that no corner of $Q$ embedds into $L(G_i),i=1,2$.
Note that $M$ decomposes in two ways as a crossed product coming from a free Bernoulli action: $$M=[*_{g_1\in G_1}(u_{g_1}L(G_2)u_{g_1}^*)]\rtimes G_1=[*_{g_2\in G_2}(u_{g_2}L(G_1)u_{g_2}^*)]\rtimes G_2.$$  
This implies, using the above discussion, that we can find $u_i\in\Cal U(M)$ and $F_i\subset G_i$ finite, for all $i\in\{1,2\}$, such that $u_iQu_i^*\subset *_{g_j\in F_j}(u_{g_j}L(G_i)u_{g_j}^*)$, whenever $\{i,j\}=\{1,2\}$.

Further we get that there exists $K_i\subset G_i$ finite and $n\geq 1$ such that $$(Q)_1\subset_{1/3}\Cal H_i=L^2(\sum_{g_{j,1},g_{j,2},..,g_{j,n}\in K_j}u_{g_{j,1}}L(G_i)u_{g_{j,2}}...L(G_i)u_{g_{j,n}}), \forall \{i,j\}=\{1,2\}.$$ 
In other words, if $P_i$ denotes the orthogonal projection onto $\Cal H_i$, then $$||x-P_i(x)||_2\leq 1/3,\forall x\in (Q)_1.$$Thus, $$||x-P_1P_2(x)||_2\leq 2/3,\forall x\in (Q)_1,$$ which leads to a contradiction as $P_1P_2$ is the orthogonal projection onto a finite dimensional Hilbert space, while $Q$ is a diffuse von Neumann algebra.

 \head 5. Central sequences in wreath products.\endhead

\vskip .1in
\noindent

In connection with Connes' $\chi(M)$ invariant, the following question has been studied in [Jo1] and [HJ]: given a crossed product $M=N\rtimes_{\sigma} G$, with $M,N$  finite von Neumann algebras, under what conditions are all central sequences of $M$ assimptotically contained in $N$? Note that if $\omega$ is a free ultrafilter on $\Bbb N$, then this is equivalent to $M'\cap M^{\omega}\subset N^{\omega}$.
If $G$ is a non-inner amenable group, then the following stronger result holds true: $L(G)'\cap M^{\omega}\subset N^{\omega}$ ([HJ]). 

In [Jo1], Jones showed that $M'\cap M^{\omega}\subset N^{\omega}$ if we assume that $N$ is a full II$_1
$ factor(i.e. has no non-trivial central sequences) and that $\varepsilon(G)\subset\text{Out}(N)$ is discrete w.r.t. the topology given by pointwise $||.||_2$ convergence. To prove this, it is first shown that $G\subset \text{Aut}(N)$ satisfies:

{\it  Condition 5.1.} There exist $ u_1,u_2..u_n\in\Cal U(N)$ and $C>0$ such that  $$\sum_{i=1}^n||\xi\sigma_g(u_i)- u_i\xi||_2^2\geq C||\xi||_2^2,\forall g\in G\setminus\{e\},\forall \xi\in L^2(N).$$
\vskip 0.05in
Below we will also consider this condition for a set (not necessarily group) of automorphisms of $N$.
In turn, (5.1.) easily implies the following:
\proclaim{ Proposition 5.2 ([Jo1])} Let $(N,\tau_N),(P,\tau_P)$ be two finite von Neumann algebras,  let $G$ be a countable, discrete group and  $\alpha:G\rightarrow \text{Aut}(N,\tau_N)$, $\beta:G\rightarrow\text{Aut}(P,\tau_P)$ be two actions. Denote by $\sigma:G\rightarrow\text{Aut}(N\overline{\otimes}P,\tau_N\times \tau_P)$ the diagonal product action and let $M=(N\overline{\otimes}P)\rtimes_{\sigma}G$. If  $\alpha(G\setminus\{e\}) \subset \text{Aut}(N,\tau)$ satisfies (5.1.),
then $N'\cap {M}^{\omega}\subset {(N\overline{\otimes}P)}^{\omega}$.
\endproclaim
{\it Proof.} Let $u_1,..,u_n\in\Cal U(N)$  and $C>0$ such that (5.1.) is satisfied. If $\{\xi_j|j\in J\}$ is an orthonormal basis for $L^2(P)$, then any $x\in M\subset L^2(M)$ can be written as $x=\sum_{j\in J,g\in G}(x_{j,g}\otimes \xi_j)u_g$, for some $x_{j,g}\in L^2(N)$. Using condition (5.1.) we get that$$\sum_{i=1}^n ||[x,u_i]||_2^2=\sum_{i=1}^n||\sum_{j\in J,g\in G}[(u_ix_{j,g}-x_{j,g}\alpha_g(u_i))\otimes\xi_j]u_g||_2^2=$$
$$\sum_{i=1}^n\sum_{j\in J,g\in G}||u_ix_{j,g}-x_{j,g}\alpha_g(u_i) ||_2^2=\sum_{j\in J,g\in G}\sum_{i=1}^n ||u_ix_{j,g}-x_{j,g}\alpha_g(u_i) ||_2^2\geq$$ $$ C\sum_{j\in J,g\in G\setminus\{e\}}||x_{j,g}||_2^2= C||x-E_{N\overline{\otimes}P}(x)||_2^2,\forall x\in M,$$
  thus the conclusion holds true.\hfill $\blacksquare$

Next, we obseve  that Bernoulli shift actions with non-amenable basis verify condition 5.1.

\proclaim{ Proposition 5.3} Let $P$ be a non-hyperfinite II$_1$-factor and denote $N=\overline{\bigotimes}_{\Bbb N}P$. For a permutation $\pi\in S(\Bbb N)$, denote by $\theta_{\pi}$ the corresponding automorphism of $N$.
 
Then for all $ n\in \Bbb N$, $S_n=\{\theta_{\pi}|\pi\in S(\Bbb N),\exists i\in\{1,2..,n\},\pi(i)\not=i\}\subset\text{Aut}(N)$
satisfies condition (5.1.).
\endproclaim 
{\it Proof.}  Recall  that by Connes' characterization of hyperfiniteness ([Co1]), a $II_1$ factor $P$ is non-hyperfinite {\it iff} there exists $u_1,u_2,..,u_n\in\Cal U(P)$ and $C>0$ such that $$\sum_{i=1}^n ||u_i\xi-\xi u_i||_2^2\geq C||\xi||_2^2,\forall \xi\in \Cal H=L^2(P)\overline{\otimes}L^2(P),$$ where $\Cal H$ is endowed with the natural $P-P$ bimodule structure, given by $a(x\otimes y)b=(ax)\otimes(yb),$ for all $a,b\in P, x,y\in L^2(P).$ 

Note that is it sufficient to prove the case $n=0$. Let $\pi\in S_0$, then $k=\pi(0)\not=0$. 
Let $\{\eta_j\}_{j\in J}\subset \overline{\bigotimes}_{\Bbb N\setminus\{0,k\}} P$ be an orthonormal basis and identify $(L^2(P))_0\overline{\otimes}(L^2(P))_k$ with $\Cal H$. Then any element $\xi\in L^2(N)$ can be written as $\xi=\Sigma_{j}x_j\otimes \eta_j$, where $x_j\in \Cal H$.
 
Thus $$\sum_{i=1}^n||\xi\theta_{\pi}(u_i)-\xi u_i||_2^2=\sum_{i=1}^n||\sum_{j\in J}(x_ju_i-u_ix_j)\otimes\eta_j||_2^2=$$
$$\sum_{i,j}||x_ju_i-u_ix_j||_2^2=\sum_{j\in J}(\sum_{i=1}^{n}||x_ju_i-u_ix_j||_2^2)\geq$$
$$C \sum_{j\in J}||x_j||_2^2=C||\xi||_2^2,\forall \pi\in S_0,\forall\xi\in L^2(N).$$ 
\hfill $\blacksquare$

If $H$ is a non-amenable group then the left-right representation of $H$ on $l^2(H\times H)$ has no almost invariant vectors, hence the previous proposition can be proved for $P=L(H)$ without making use of Connes' powerful result. 
\proclaim{Corollary 5.4} Let $P$ be a non-hyperfinite $II_1$ factor, $G$ be a countable discrete group and denote $M=P\wr G$. Then $P'\cap M^{\omega}\subset (\overline{\bigotimes}_{G}P)^{\omega}$, thus,  $M'\cap M^{\omega}\subset (\overline{\bigotimes}_{G}P)^{\omega}$.
\endproclaim

\head 6. Constructing central sequences via Ornstein-Weiss\endhead

In the previous section we have shown that if $N$ is a non-amenable II$_1$ factor then for every countable discrete group $G$, central sequences in $N\wr G$ have to assimptoticaly lie in $\overline{\bigotimes}_{G}N$. On the other hand, if $G$ is amenable, then  we  can
use Ornstein-Weiss' Rokhlin's lemma to produce many sequences in $\overline{\bigotimes}_{G}N$, which are central in $N\wr G$. Using this, we will able to recover the core $\overline{\bigotimes}_{G}N$.

Recall  the following consequence of  Ornstein-Weiss' Rokhlin lemma: 

\proclaim{Theorem 6.1.([OW])} Let $G$ be an amenable group acting on standard probability space $(X,\mu)$ by measure preserving transformations. Then for all $\varepsilon>0$ and $F\subset G$ finite, we can find  $K_1,..,K_n$ finite $(F,\varepsilon)-$invaried subsets of $G$ (i.e. $|K_i\Delta gK_i|\leq\varepsilon|K_i|,$ for all $g\in F,i=1,..,n$) and measurable sets $B_1,..,B_n\subset X$ such that 
 
$(i)\{K_iB_i|1\leq i\leq n\}$ are disjoint,

$(ii)\{tB_i|t\in K_i\}$ are disjoint, $1\leq i\leq n$.

$(iii)\mu(\cup K_iB_i)>1-\varepsilon$.
\endproclaim
\vskip 0.1in

\proclaim {Theorem 6.2} Let $N$ be a II$_1$ factor, $G$ be an infinite amenable group and denote $M=N\wr G$. Let $\tilde M\supset M$ be a II$_1$ factor such that $M$ and $(\overline{\bigotimes}_{G}N)'\cap \tilde M$ generate $\tilde M$ (e.g. if $\tilde M=M$). Assume that $Q\subset \tilde M$ is a regular von Neumann subalgebra     such that $\tilde M'\cap {\tilde M}^{\omega}\subset Q^{\omega}$. Then a corner of $\overline{\bigotimes}_{G}N$ can be embedded into $Q$ inside $\tilde M$.
\endproclaim 
{\it Proof.} Using the hypothesis we deduce that for all $ \varepsilon>0$, there exists $\delta=\delta(\varepsilon)>0$ and $F=F(\varepsilon)\subset G$ finite such that if $U\in (\tilde M)_1$ satisfies $$U\in_{\delta}\overline{\bigotimes}_{G\setminus F}N,||\sigma_g(U)-U||_2\leq\delta,\forall g\in F,$$ then $$||U-E_Q(U)||_2\leq \varepsilon.$$
\vskip 0.05in
 
 Fix $1>\varepsilon>0$ and let $\delta,F$ be defined as above.
If $A\subset N$ is a MASA,  then $G$ acts on $\overline{\bigotimes}_{G}A$. By applying the previous theorem to this action we obtain that there exists projections $p_i\in\Cal P(\overline{\bigotimes}_{G}A)$ and $(F,{\delta}^2/4)$-invaried subsets $K_1,..,K_n$ of $G$ such that the projections $\{\sigma_{k_i}(p_i)\}_{i=\overline{1,n},k_i\in K_i}$ are mutually disjoint and satisfy $$
1-\varepsilon+{\varepsilon}^2\leq\tau(\sum_{i=1}^n\sum_{k_i\in K_i} \sigma_{k_i}(p_i))=\sum_{i=1}^{n}|K_i|\tau(q_i)\leq 1.$$
 
 Next, we can find a finite set $K\subset G$ such that $p_i\in_{\delta \tau(p_i)/2}\overline{\bigotimes}_K A,$ for all $i=1,n$. Let $g\in G\setminus (K^{-1}(\cup K_i)^{-1}F)$ and denote by $\rho$ the right Bernoulli shift action of $G$ on $\overline{\bigotimes}_G N$. If we set $q_i=\rho_g(p_i)$, then $\sigma_{k_i}(q_i)\in_{\delta\tau(q_i)}\overline{\bigotimes}_{K_iKh} N\subset \overline{\bigotimes}_{G\setminus F}N$ and  $\{\sigma_{k_i}(q_i)\}_{i=\overline{1,n},k_i\in K_i}$ is  an $(\varepsilon-{\varepsilon}^2)$- almost partition of unity, in the same sense as above. 
  
Now, let $u\in\Cal U(\overline{\bigotimes}_{G\setminus [(\cup K_i)^{-1}F]}N)$ and set $U=\sum_{i=1}^n\sum_{k_i\in K_i}\sigma_{k_i}(q_iuq_i)$.
Using the fact that $|K_i\Delta gK_i|/|K_i|\leq {\delta^2}/4,$ for all $i=\overline{1,n},$ and all $g\in F$, we have the following estimates:
$$||\sigma_g(U)-U||_2^2\leq 4||\sum_{i=1}^n\sum_{k_i\in K_i\setminus gK_i}\sigma_{k_i}(q_iuq_i)||_2^2\leq$$ 
$$ 4\sum_{i=1}^n |K_i\setminus gK_i|\tau(q_i)\leq 4\max_{i=\overline{1,n}}(|K_i\setminus gK_i|/|K_i|)\sum_{i=1}^n |K_i|\tau(q_i)\leq$$ $$ 4({\delta}^{2}/4)={\delta}^2,\forall g\in F.$$
 Since $\sigma_{k_i}(q_i)\in_{\delta\tau(q_i)/2}\overline{\bigotimes}_{G\setminus F}N,$ for all $i=\overline{1,n},$ and all $k_i\in K_i$, we get  that $$U\in_{\sum_{i=1}^n\sum_{k_i\in K_i}\delta\tau(q_i)} \overline{\bigotimes}_{G\setminus F}N,$$
hence that $U\in_{\delta}\overline{\bigotimes}_{G\setminus F}N$.  Thus, we must have that $||E_Q(U)-U||_2\leq\varepsilon$,or , equivalently, that 
\vskip 0.05in
{\it Fact 1.} If  $V=(\cup K_i)^{-1}F$, then $\forall u\in\Cal U(\overline{\bigotimes}_{G\setminus V}N)$ we have that $$||E_{Q}(\sum_{i=1}^{n}\sum_{k_i\in K_i} \sigma_{k_i}(q_{i}uq_i))-\sum_{i=1}^{n}\sum_{k_i\in K_i} \sigma_{k_i}(q_iuq_i)||_2\leq \varepsilon.$$ 
\vskip 0.05in
Since $\lim_{W\nearrow G}\sup_{u\in\Cal U(\overline{\bigotimes}_{G\setminus W}N)}||[u,q]||_2=0,$ for all $q\in \overline{\bigotimes}_G N,$ we deduce that	there exists $W\supset V$ finite such that $$||\sum_{i=1}^{n}\sum_{k_i\in K_i} \sigma_{k_i}(q_iuq_i)||_2^2\geq \sum_{i=1}^n\sum_{k_i\in K_i}\tau(q_i)-\varepsilon\geq (1-\varepsilon)^2,\forall u\in\Cal U(\overline{\bigotimes}_{G\setminus W}N).$$ Let $\{v_i|i=\overline{1,m}\}=\{U_{k_i}q_i|i=\overline{1,n},k_i\in K_i\}$, where we denote by $\{U_g|g\in G\}\subset \Cal U(M)$ the canonical unitaries  implementing the action of $G$ on $\overline{\bigotimes}_{G}N$. Then the last inequality implies: for all $\varepsilon>0$, there exists $W(\varepsilon)\subset G$ finite, $m=m(\varepsilon)\in\Bbb N$ and $\{v_i\}_{i=1,m}$ partial isometries with mutually orthogonal supports such that $$||E_Q(\sum_{i=1}^m v_iu {v_i}^*)||_2 \geq 1-2\varepsilon, \forall u\in\Cal U(\overline{\bigotimes}_{G\setminus W(\varepsilon)}N).$$

\vskip 0.05in
Let $\{\eta_l\}_{l\geq 0}$ be an orthonormal basis of $M$ over $Q$ and for all $l\geq 0$, denote by  $P_l$   the orthogonal projection onto $L^2(\sum_{i=1}^{l}Q\eta_i)$.
\vskip 0.05in 

{\it Fact 2.} For all $\varepsilon>0$, there exists $l(\varepsilon)$ such that $$\overline{\bigotimes}_{G\setminus W(\varepsilon)}N\subset_{3\varepsilon}\sum_{i=0}^{l(\varepsilon)}Q\eta_i.$$
\vskip 0.05in
  To prove this, first note that since $Q\subset \tilde M$ is regular, we can find $l(\varepsilon)$  such that $$||E_Q(\sum_{i=1}^m v_iu{v_i}^*)-E_{Q}(\sum_{i=1}^m v_i P_{l(\varepsilon)}(u){v_i}^*)||_2\leq\varepsilon ,\forall u\in\Cal U(M).$$ By combining this with the above, we get that $$||E_{Q}(\sum_{i=1}^m v_i P_{l(\varepsilon)}(u){v_i}^*)||_2\geq 1-3\varepsilon,\forall u\in\Cal U(\overline{\bigotimes}_{G\setminus W(\varepsilon)}N).$$
On the other hand, we have $$ ||E_{Q}(\sum_{i=1}^m v_i P_{l(\varepsilon)}(u){v_i}^*)||_2\leq ||\sum_{i=1}^m v_i P_{l(\varepsilon)}(u){v_i}^*||_2\leq ||P_{l(\varepsilon)}(u)||_2,\forall u\in M,$$ which entails $$||P_{l(\varepsilon)}(u)||_2\geq 1-3\varepsilon,\forall u\in\Cal U(\overline{\bigotimes}_{G\setminus W(\varepsilon)}N).$$

\vskip 0.05in
{\it  Fact 3.} For all $\varepsilon>0$, there exists $k(\varepsilon)\geq l(\varepsilon)$ such that $$\overline{\bigotimes}_{W(\varepsilon)}N\subset_{4\varepsilon}\sum_{n=1}^{k(\varepsilon)}Q\eta_i.$$
\vskip 0.05in

This fact is trivial since we can find $g\in G$ such that $W(\varepsilon)g^{-1}\subset G\setminus W(\varepsilon)$ and thus, by fact 2, $$ \overline{\bigotimes}_{W(\varepsilon)}N\subset_{3\varepsilon}(\sum_{i=1}^{l(\varepsilon)}Q\eta_i)U_g.$$
\vskip 0.05in

Finally, note that by the last two facts $$(\overline{\bigotimes}_{G\setminus W(\varepsilon)}N),(\overline{\bigotimes}_{W(\varepsilon)}N)\subset_{4\varepsilon}\sum_{i=1}^{k(\varepsilon)}Q\eta_i,\forall \varepsilon>0.$$
If $\varepsilon<1/8$, then $4\varepsilon+4\varepsilon<1$, hence   we can apply Corollary 1.2.(ii) to deduce that a corner of $\overline{\bigotimes}_{G}N$ embedds into $Q$.
\hfill$\blacksquare$

\proclaim{Theorem 6.3} 
Let $G,H$ be two ICC   groups, one of them being amenable, let $P$ and $Q$ be two non-hyperfinite II$_1$ factors and denote $M=P\wr G$, $N=Q\wr H$.  
If $\theta:M\rightarrow N$ is an isomorphism, then there exists $u\in\Cal U(N)$ such that $u\theta(\overline{\bigotimes}_{G}P)u^*=\overline{\bigotimes}_{H}Q$. In particular, $G\simeq H$.
\endproclaim
{\it Proof}. After identification via $\theta$ we may assume that $M=N$. Since $M$ (respectively $N$) has non-trivial central sequences iff $G$ (respectively $H$) is amenable, we deduce that both groups are amenable.
Since $P$ and $Q$ are  non-hyperfinite, we can apply Corollary 2.4. to deduce that   $M'\cap M^{\omega}\subset {(\overline{\bigotimes}_{G}P)}^{\omega}$ and that $M'\cap M^{\omega}\subset {(\overline{\bigotimes}_{H}Q)}^{\omega}$. 

Using these inclusions and the fact that both $\overline{\bigotimes}_{G}P$ and $\overline{\bigotimes}_{H}Q$ are regular subfactors of $M$, we can apply Theorem 6.2. to get  two $\overline{\bigotimes}_{G}P-\overline{\bigotimes}_{H}Q$ Hilbert bimodules $\Cal H,\Cal K\subset L^2(M)$, which are finite dimensional over $\overline{\bigotimes}_{G}P$ and over $\overline{\bigotimes}_{H}Q$, respectively. 

Since $G$ and $H$ are ICC groups, we can use Theorem 8.4. in [IPP]  to get the conclusion.
\hfill $\blacksquare$

\vskip 0.1in
For a finite factor $M$, we denote by $\sigma_{M}$ the Bernoulli shift action of $G$ on $\overline{\bigotimes}_{G}M$; for a cocycle action $\sigma$ of $G$ on a II$_1$ factor $M$ and for $t>0$, we denote by $\sigma^t$ the induced cocycle action of $G$ on $M^t$. If two cocycle actions $\alpha$ and $\beta$ are cocycle conjugate, then we write $\alpha\sim\beta$.  Also, for two actions $\alpha$,$\beta$ of the same group $G$ on $M$ and $N$, respectively, we denote by $\alpha\otimes\beta$ the diagonal product action of $G$ on $M\overline{\otimes}N$. 
\proclaim{Proposition 6.4} Let $N$ be a type II$_1$ factor and let $G$ be an $\infty$ amenable group. Then 

(i) $\Cal F(N\wr G)=\Bbb R_{+}^*$.

(ii) $(\overline{\bigotimes}_{G}N\subset N\wr G)\simeq (\overline{\bigotimes}_{G} (N\overline{\otimes}R)\subset (N\overline{\otimes}R)\wr G),$ for all $t>0$.
\endproclaim
{\it Proof.} 
Recall that by a result of Ocneanu ([Oc]) any two cocycle actions of an amenable $G$ group on the hyperfinite II$_1$ factor $R$ are cocycle conjugate. Thus, if $\alpha,\beta:G\rightarrow\text{Aut}(R)$ are cocycle actions and $\sigma:G\rightarrow\text{Aut}(M)$ is an action on  a von Neumann algebra $M$, then the cocycle actions $\alpha\otimes\sigma$ and $\beta\otimes\sigma$ on $R\overline\otimes M$
are cocycle conjugate. Also, note that if $\alpha,\beta$ are two cocycle actions of the same group $G$, then $(\alpha\otimes\beta)^t\sim{\alpha}^t\otimes \beta$.

If $N$ is a II$_1$ factor, then by combining the above statements we get that
 $\sigma_{N^t}=(\sigma_{N^{1/2}}\otimes\sigma_{\Bbb M_2(\Bbb C)})^t\sim (\sigma_{N^{1/2}}\otimes \sigma_{R})^t\sim\sigma_{N^{1/2}}\otimes\sigma_{R}^t\sim\sigma_{N^{1/2}}\otimes\sigma_{R}\sim\sigma_N,$ hence, in particular, $\Cal F(N\wr G)=\Bbb R_{+}^*$. The second claim follows in a similar way: $\sigma_{N}\sim \sigma_{N^{1/2}}\otimes\sigma_{\Bbb M_2(\Bbb C)}\sim\sigma_{N^{1/2}}\otimes\sigma_{R}\sim\sigma_{N^{1/2}\overline{\otimes}R}\sim\sigma_{N\overline{\otimes}R}$.    \hfill $\blacksquare$
\vskip 0.05in
Since a result of Dykema-Radulescu gives that for all $m,n\geq 2,$ there exists $t>0$ such that $L(\Bbb F_m)^t\simeq 
L(\Bbb F_n)$, by using the previous two results we get the following:
\vskip 0.05in

\proclaim{Corollary 6.5} If $2\leq m,n<\infty$ and $G,H$ are two $\infty$  ICC amenable groups, then $L(\Bbb F_m\wr G)\simeq L(F_n\wr H)$ {\it iff} $G\simeq H$.
\endproclaim


\head  References\endhead

\item {[Co1]} A.Connes: {\it Classification of injective factors}, Ann. of Math., {\bf 104} (1976), 73-115.
\item {[Co2]} A. Connes: {\it Sur la classification des facteurs de type \rm II}, C. R. Acad. Sci. Paris SŽr.I Math.
{\bf 281} (1975), A13ÐA15.

\item {[HJ]} R.H. Herman, V.F.R Jones, {\it Central sequences in crossed products. }Operator algebras and mathematical physics (Iowa City, Iowa, 1985), 539--544, 
Contemp. Math., 62, Amer. Math. Soc., Providence, RI, 1987.
\item {[Jo1]} V.F.R. Jones: {\it Central sequences in crossed products of full factors}, Duke Math. J.  {\bf 49}, no. 1 (1982), 29–33 
\item {[Jo2]} V.F.R. Jones: {\it Notes on Connes' invariant $\chi(M)$}, unpublished.
\item {[IPP]} A. Ioana, J. Peterson, S. Popa: {\it Amalgamated Free Products of w-Rigid Factors and Calculation of their Symmetry Groups}, math.OA/0505589 , Preprint 2005. 
\item {[MV]} F. Martin, A. Valette: {\it On the first Lp-cohomology of discrete groups}
\item {[Oc]} A. Ocneanu: {\it Actions of discrete amenable groups on von Neumann algebras}, Lecture Notes in Mathematics, 1138, Springer Verlag, Berlin, 1985. 
\item {[OW]} D. Ornstein, B.Weiss: {\it Ergodic theory of amenable groups. I. The Rokhlin lemma.}, Bull. Amer. Math. Soc. (N.S.) {\bf 1} (1980), 161-164.

\item {[OP]} N. Ozawa, S. Popa: {\it Some prime factorization results for type \rm II$\sb 1$ factors}, Invent. Math., {\bf 156} (2004), 223-234. 
\item {[Po1]} S. Popa: {\it On a class of type II$_1$ factors with Betti numbers invariants}, math. OA/0209130, to appear in Annals of Math. 
\item {[Po2]} S. Popa: {\it Strong rigidity of II$_1$ factors arising from malleable actions of w-rigid groups}, Part I, math.OA/0305306, Preprint 2003, to appear in Invent. Math.

\item {[Po3]}S. Popa: {\it Strong rigidity of II$_1$ factors arising from malleable actions of w-rigid groups},
Part II, math.OA/0407103, Preprint 2004, to appear in Invent. Math.

\item {[Po4]} S. Popa: {\it A unique decomposition result for HT factors with torsion free core}, math.OA/0401138, Preprint 2004.

\item {[Po5]} S. Popa: {\it Some rigidity results for non-commutative Bernoulii shifts}, Journal of Functional Analysis {\bf 230} (2005), 1-56.

\item {[Po6]} S. Popa: {\it Cocycle and orbit equivalence superrigidity for Bernoulli actions of Kazhdan groups} 
 
\item {[Va]} S. Vaes: {\it Rigidity results for Bernoulli actions and their von Neumann algebras}

\enddocument